\documentclass[reqno, openany, amssymb, psamsfonts]{amsart}
\usepackage{amsrefs}
\usepackage{amsmath}
\usepackage{amsthm}
\usepackage{mathrsfs,comment}
\usepackage[usenames,dvipsnames]{color}
\usepackage[normalem]{ulem}
\usepackage{url}
\usepackage[all,arc,2cell]{xy}
\UseAllTwocells
\usepackage{enumerate}
\usepackage{hyperref}
\usepackage{tabularx}
\usepackage[pdftex]{graphicx}
\usepackage{float}

\calclayout

\calclayout

\hypersetup{%
  bookmarksnumbered=true,%
  bookmarks=true,%
  colorlinks=true,%
  linkcolor=blue,%
  citecolor=blue,%
  filecolor=blue,%
  menucolor=blue,%
  pagecolor=blue,%
  urlcolor=blue,%
  pdfnewwindow=true,%
  pdfstartview=FitBH}

\def\makeautorefname#1#2{\expandafter\def\csname#1autorefname\endcsname{#2}}

\def\equationautorefname~#1\null{(#1)\null}
\makeautorefname{footnote}{footnote}%
\makeautorefname{item}{item}%
\makeautorefname{figure}{Figure}%
\makeautorefname{table}{Table}%
\makeautorefname{part}{Part}%
\makeautorefname{appendix}{Appendix}%
\makeautorefname{chapter}{Chapter}%
\makeautorefname{section}{Section}%
\makeautorefname{subsection}{Section}%
\makeautorefname{subsubsection}{Section}%
\makeautorefname{theorem}{Theorem}%
\makeautorefname{thm}{Theorem}%
\makeautorefname{cor}{Corollary}%
\makeautorefname{lem}{Lemma}%
\makeautorefname{prop}{Proposition}%
\makeautorefname{pro}{Property}
\makeautorefname{conj}{Conjecture}%
\makeautorefname{defn}{Definition}%
\makeautorefname{notn}{Notation}
\makeautorefname{notns}{Notations}
\makeautorefname{rem}{Remark}%
\makeautorefname{quest}{Question}%
\makeautorefname{exmp}{Example}%
\makeautorefname{ax}{Axiom}%
\makeautorefname{claim}{Claim}%
\makeautorefname{ass}{Assumption}%
\makeautorefname{asss}{Assumptions}%
\makeautorefname{con}{Construction}%
\makeautorefname{prob}{Problem}%
\makeautorefname{warn}{Warning}%
\makeautorefname{obs}{Observation}%
\makeautorefname{conv}{Convention}%

\newtheorem{thm}{Theorem}[section]

\newtheorem{lem}{Lemma}[section]

\theoremstyle{definition}
\newtheorem{defn}{Definition}[section]
\newtheorem{ass}{Assumption}[section]

\newtheorem{pro}{Property}[section]

\newtheorem{rem}{Remark}[section]

\makeatletter
\let\c@obs=\c@thm
\let\c@cor=\c@thm
\let\c@prop=\c@thm
\let\c@lem=\c@thm
\let\c@prob=\c@thm
\let\c@con=\c@thm
\let\c@conj=\c@thm
\let\c@defn=\c@thm
\let\c@notn=\c@thm
\let\c@notns=\c@thm
\let\c@exmp=\c@thm
\let\c@ax=\c@thm
\let\c@pro=\c@thm
\let\c@ass=\c@thm
\let\c@warn=\c@thm
\let\c@rem=\c@thm
\let\c@sch=\c@thm
\numberwithin{equation}{section}
\makeatother

\newcommand{\ZZ}{\mathbb{Z}}
\newcommand{\RR}{\mathbb{R}}

\usepackage{fancyhdr}
\def\roq#1{\ \quad \textnormal{#1}\quad \ }
\def\O{{\rm O}\mskip1mu}
\def\:{\mskip2mu}

\newcommand{\iu}{\mathrm{i}}

\DeclareMathOperator{\sgn}{sgn}

\title{Continuous Revival of the Periodic Schr\"odinger Equation\break with Piecewise $C^2$ Potentials}

\author{Dinh-Quan Tran$^1$}
\address{$^1$School of Mathematics, University of Minnesota, Minneapolis, MN 55455}
\email{dinhquan.tran23@gmail.com}

\author{Peter J. Olver$^2$}
\address{$^2$School of Mathematics, University of Minnesota, Minneapolis, MN 55455}
\email{olver@umn.edu}
\urladdr{http://www.math.umn.edu/$\sim$olver}

\begin{document}

\begin{abstract}

In this paper, we investigate the revivals of the one-dimensional periodic Schr\"odinger equation with a piecewise $C^2$ potential function. As has been observed through numerical simulations of the equation with various initial data and potential functions, the solution, while remaining fractalized at irrational times, exhibits a form of revival at rational times. The goal is to prove that the solution at these rational times is given by a finite linear combination of translations and dilations of the initial datum, plus an additional continuous term, which we call ``continuous revival". In pursuit of this result, we present a review of relevant properties of the periodic Schr\"odinger equation as an eigenvalue problem, including asymptotic results on both the eigenvalues and eigenfunctions.

\end{abstract}

\maketitle

\vspace{-3ex}


\setlength{\parskip}{3pt}

\vspace{3ex}

\section{Introduction} 

The Schr\"odinger equation, a partial differential equation that governs the wave function of non-relativistic quantum-mechanical systems, has been a cornerstone of mathematical physics. Named after the twentieth century physicist Erwin Schr\"odinger, the equation formed the basis of his Nobel Prize and continues to be of fundamental importance for applications, particularly in quantum mechanics. Within this paper, we are mainly concerned with the one-dimensional periodic Schr\"odinger equation.

\begin{defn}[Periodic Schr\"odinger equation] \label{defn:periodic_schrodinger}
Consider an initial datum $f \in L^2[0, 2\pi]$ and a real-valued potential function $V(x) \in L^\infty[0, 2\pi]$. A solution $u(t, x)$ of the periodic Schr\"odinger equation satisfies the initial value problem 
\begin{align} \label{eq:schrodinger}
    \begin{cases}
        \iu \, u_t = - u_{xx} + Vu, \\[1pt]
        u(0, x) = f(x),
    \end{cases}
\end{align}
where $\iu = \sqrt{-1}$, subject to periodic boundary conditions
\begin{align} \label{eq:periodic_boundary}
    \begin{cases}
        u(t, 0) = u(t, 2\:\pi), \\[1pt]
        u_x(t, 0) = u_x(t, 2\:\pi).
    \end{cases}
\end{align}
Here, with
\begin{equation} \label{eq:schrodinger_operator}
    L := -\frac{d^2}{dx^2} + V
\end{equation}
being the Schr\"odinger operator, we take the mild solution 
\[ u(x) = e^{-itL} f(x) \]
as an element of $L^2[0, 2\pi]$, and all identities involving $u$ in this paper shall henceforth be regarded in an almost everywhere sense.
\end{defn}

As a dispersive equation on a periodic domain, the Schr\"odinger equation displays two vastly different behaviors depending on whether the time is a rational multiple of the period. While the solution remains fractalized at irrational times, it exhibits revival at rational times, being related to a finite linear combination of translations and dilations of the initial data. This quantization phenomenon was first observed by Henry Fox Talbot in 1836 when he noticed that the image of a diffraction grating is repeated at regular distances away from the grating plane. Since then, the phenomenon, called the ``Talbot effect" \cite{berry_1996, berry_marzoli_schleich_2001}, or ``dispersive quantization" \cite{olver_2010}, has been noted in various periodic dispersive equations, including the linearized Korteweg--deVries equation (also known as the Airy equation) \cite{olver_2010}, the free space Schr\"odinger equation without potential \cite{berry_1996,chen_olver_2014}, and related linear and nonlinear equations \cite{chen_olver_2014, chousionis_erdogan_tzirakis_2015, erdogan_tzirakis_2016}. More recently, the continuous revival phenomenon has been investigated for the class of Schr\"odinger equations with $H^2$ potential \cite{boulton_farmakis_pelloni_2024}, while the fractal dimension of the solution to the equation at irrational times has also been analyzed for a broader class of $H^s$ potentials \cite{cho_2024}. In order to extend our understanding to more general potential functions, it is important that we first characterize revival for potential functions with discontinuities. This paper, therefore, aims to investigate the phenomenon in a simpler setting, where the potential $V(x)$ is piecewise $C^2$.

\begin{ass}[Piecewise $C^2$ potential function] \label{ass:potential}
The potential function $V(x)$ is piecewise $C^2$ if there exist finitely many points $0 = x_0 < x_1 < x_2 < ... < x_N < x_{N+1} = 2\:\pi$ such that 
\begin{itemize}
    \item $V(x)$ is bounded and $2\:\pi$-periodic;
    \item $V(x) \in C^2(x_i, x_{i+1})$ for all $0 \leq i \leq N$;
    \item $V(x), V'(x)$, and $V''(x)$ have finite right and left limits at the discontinuities $x_i$. 
\end{itemize}
\end{ass}

Note that since changing the potential $V$ at finitely many points does not change the associated Schr\"odinger operator, potentials that agree almost everywhere are identified. Accordingly, inequalities between potentials are understood to hold almost everywhere unless otherwise stated. 

The central focus of our project shall be the following theorem mirroring that of Boulton et al. \cite{boulton_farmakis_pelloni_2024}.

\begin{thm}[Continuous revival at rational times] \label{thm:revival}
Consider an initial datum $f \in L^2[0, 2\pi]$, and let $f^*$ be the $2\:\pi$-periodic extension of $f$. Furthermore, let $u(t, x)$ be the solution of the periodic Schr\"odinger equation \eqref{eq:schrodinger} with the given initial datum $f$ and a potential function $V(x)$ satisfying Assumption \ref{ass:potential}. For $q, r \in \mathbb{N}$ coprime numbers, define the revival function as
\begin{equation} \label{eq:revival}
    \psi \left( 2\:\pi \frac{q}{r}, x \right) = \frac{1}{r} \,e^{-2\:\pi \iu \langle V \rangle q/r} \sum^{r-1}_{k, m=0} e^{2\:\pi \iu  \left(mk/r- m^2q/r \right)} f^*\left(x - 2\:\pi \frac{k}{r} \right).
\end{equation}
Then, at rational times $t_0 = 2\:\pi \frac{q}{r}$, there exists a function $w(t_0,x)$ continuous in $x$ such that in $L^2[0, 2\pi]$,
\begin{equation*}
    u \left( t_0, x \right) = w \left( t_0, x \right) + \psi \left( t_0, x \right).
\end{equation*}
\end{thm}

\section{Sturm-Liouville theory for the Schr\"odinger equation}

Inserting the ansatz $u(t, x) = e^{-\iu \lambda t} \psi(x)$ into \eqref{eq:schrodinger} produces the ordinary differential equation 
\begin{align} \label{eq:eigenvalue_relation}
    \lambda \psi = - \,\psi'' + V \psi = L[\psi]
\end{align}
along with the periodic boundary conditions
\begin{align} \label{eq:eigenfunc_periodic_boundary_condition}
    \begin{cases}
        \psi(0) = \psi(2\:\pi), \\[1pt]
        \psi'(0) = \psi'(2\:\pi).
    \end{cases}
\end{align}
Here, the Schr\"odinger operator $L$ as in \eqref{eq:schrodinger_operator} is an operator on $L^2[0, 2\pi]$ with domain
\[ D(L) = \{ \psi \in H^2[0, 2\pi]: \psi(0) = \psi(2\:\pi), \psi'(0) = \psi'(2\:\pi)\}. \]
The operator $L$ is self-adjoint with compact resolvent, and hence the eigenvalue equation \eqref{eq:eigenvalue_relation} admits a countably infinite sequence of real eigenvalues
\[ \lambda_0 \leq \lambda_1 \leq \lambda_2 \leq ... \roq{with } \lambda_n \to \infty \ \roq{as} \ n \to \infty, \]
with multiplicities included. Using these eigenvalues, we can choose a basis of eigenfunctions to be real-valued and orthonormal on $[0, 2\:\pi]$; that is
\begin{align*}
    \int_0^{2\:\pi} \psi_m(x) \psi_n(x) \, dx = \begin{cases} 
        1 \qquad &\roq{if} m = n,\\
        0 &\roq{if} m \neq n.
    \end{cases}
\end{align*}
These eigenfunctions can be trivially extended into continuously differentiable functions on the entirety of the real line thanks to the periodic boundary conditions \eqref{eq:eigenfunc_periodic_boundary_condition}. Thus, for an initial datum $f \in L^2[0, 2\pi]$, the unique solution $u \in L^2[0, 2\pi]$ is given by
\[ u(t, \cdot) = e^{-itL} f = \sum_{n = 0}^\infty \langle f, \psi_n \rangle e^{-i \lambda_n t} \psi_n. \]
Our goal for this section is to quantify the asymptotic behavior of the eigenvalues and eigenfunctions of \eqref{eq:eigenvalue_relation}, which shall be important for our proof of Theorem \ref{thm:revival}. To prepare for those results, we shall, in the next subsections, look at some variational results and comparison theorems that shall give us insights into the properties of the eigenfunctions.

\subsection{Variational results}

Consider the Dirichlet integral given by
\begin{align} \label{eq:dirichlet_integral}
    J(f, g) = \int_0^{2\:\pi} f'(x) g'(x) + V(x) f(x) g(x) \, dx.
\end{align}
Applying integration by parts and the boundary conditions \eqref{eq:eigenfunc_periodic_boundary_condition}, we obtain
\begin{align*}
    J(f, g) &= - \int_0^{2\:\pi} f(x) \left[ g''(x) - V(x) g(x) \right] dx + f(x) g'(x) \bigg |_{x=0}^{x = 2\:\pi} \\
    &= - \int_0^{2\:\pi} f(x) \left[ g''(x) - V(x) g(x) \right] dx
\end{align*}
where the boundary terms vanish due to the periodicity of $f$ and $g$. Notably, we have
\begin{align} \label{eq:dirichlet_integral_eigenfunc}
    J(f, \psi_n) &= - \int_0^{2\:\pi} f(x) \left[ \psi_n''(x) - V(x) \psi_n(x) \right] dx = \int_0^{2\:\pi} f(x) \left[ \lambda_n \psi_n(x) \right] dx = \lambda_n f_n
\end{align}
where $f_n$ is the $n$-th eigenfunction expansion coefficient of $f$. An immediate corollary is
\begin{align} \label{eq:dirichlet_integral_eigenfunc_orthonormal}
    J(\psi_m, \psi_n) = \begin{cases} 
        \lambda_n \qquad &\roq{if} m = n, \\
        0 &\roq{if} m \neq n.
    \end{cases}
\end{align}

The following result gives an inequality involving $J$ that elaborates on the familiar Parseval identity.

\begin{lem} \label{lem:dirichlet_integral_parseval}
Let $f$ be a real-valued function satisfying the boundary conditions \eqref{eq:eigenfunc_periodic_boundary_condition} that is continuous with piecewise continuous derivative in $[0, 2\:\pi]$. Then, with the Fourier coefficients $f_n$ defined above,
\[ \sum_{n=0}^\infty \lambda_n f_n^2 \leq J(f, f). \]
\end{lem}

\begin{proof}
We first suppose that $V(x) \geq 0$ almost everywhere in $[0, 2\pi]$. Then, for any continuous $g(x)$ with piecewise continuous derivative, we have
\[ J(g, g) = \int_0^{2\:\pi} [g'(x)]^2 + V(x) g^2(x) \, dx \geq 0. \]
In particular, for any positive integer $N$, we have
\begin{align*}
    0 \leq J \left(f - \sum_{n=0}^N f_n \psi_n, f - \sum_{n=0}^N f_n \psi_n \right) &= J(f, f) - 2 \sum_{n=0}^N f_n J(f, \psi_n) + J \left( \sum_{n=0}^N f_n \psi_n, \sum_{n=0}^N f_n \psi_n \right) \\
    &= J(f, f) - 2 \sum_{n=0}^N \lambda_n f_n^2 + \sum_{n=0}^N \lambda_n f_n^2 \\
    &= J(f, f) - \sum_{n=0}^N \lambda_n f_n^2
\end{align*}
where, in the second-to-last equality, we used \eqref{eq:dirichlet_integral_eigenfunc} and \eqref{eq:dirichlet_integral_eigenfunc_orthonormal}. This gives the desired inequality for the case where $V(x)$ is positive, and we now aim for the general case. Let $v_0$ be a positive constant large enough such that 
\[ V(x) \geq -v_0 \roq{for a.e.} x \in [0, 2\pi] \]
Now, the original relation \eqref{eq:eigenvalue_relation} can be rewritten as
\[ \psi''(x) + (\Lambda - \widetilde{V}(x)) \psi(x) = 0 \]
where $\Lambda = \lambda + v_0$ and $\widetilde{V}(x) = V(x) + v_0$ is a positive function. Applying the first part of the proof gives
\[ \sum_{n=0}^\infty (\lambda_n + v_0) f_n^2 \leq \int_0^{2\:\pi} [f'(x)]^2 + (V(x) + v_0) f^2(x) \, dx = J(f, f) + v_0 \int_0^{2\:\pi} f^2(x) \, dx. \]
However, thanks to Parseval's formula,
\[ \sum_{n=0}^\infty f_n^2 = \int_0^{2\:\pi} f^2(x) \, dx, \]
we can cancel the extra terms related to $v_0$. The general case follows as a result.
\end{proof}

Lemma \ref{lem:dirichlet_integral_parseval} serves to establish the next variational result, which determines the changes in the eigenvalues with respect to that of the potential function.

\begin{thm} \label{thm:eigenvalue_variation}
Let 
\[ V_1(x) \geq V(x) \roq{for a.e.} x \in [0, 2\pi] \]
and let $\lambda_{1, n}$ denote the eigenvalues of the periodic problem \eqref{eq:eigenvalue_relation} with potential $V_1(x)$.
Then
\[ \lambda_{1, n} \geq \lambda_n \roq{for all} n \geq 0. \]
\end{thm}

\begin{proof}
Let $\psi_{1, n}(x)$ denote the eigenfunction in the orthonormal basis of the modified problem that corresponds to $\lambda_{1, n}$, and let $J_1(f, g)$ denote the Dirichlet integral \eqref{eq:dirichlet_integral} but with $V(x)$ replaced by $V_1(x)$. We obtain the following relation between the Dirichlet integrals:
\[ J_1(f, f) = \int_0^{2\:\pi} [f'(x)]^2 + V_1(x) f^2(x) \, dx \geq \int_0^{2\:\pi} [f'(x)]^2 + V(x) f^2(x) \, dx = J(f, f). \]
First, we seek to prove the theorem for $n = 0$. Since $\lambda_0$ is the smallest of the eigenvalues, applying Lemma \ref{lem:dirichlet_integral_parseval} to any function $f$ gives us
\[ J(f, f) \geq \lambda_0 \sum_{n=0}^\infty f_n^2 = \lambda_0 \int_0^{2\:\pi} f^2(x) \, dx. \]
Consequently,
\[ \lambda_{1, 0} = J_1(\psi_{1, 0}, \psi_{1, 0}) \geq J(\psi_{1, 0}, \psi_{1, 0}) \geq \lambda_0 \int_0^{2\:\pi} \psi_{1, 0}^2 \, dx = \lambda_0 \]
which completes the case of $n = 0$. For $n = 1$, we consider a function
\[ f(x) = c_0 \psi_{1, 0}(x) + c_1 \psi_{1, 1}(x) \]
where the real constants $c_0$ and $c_1$ are such that
\[ \begin{cases}
    c_0^2 + c_1^2 = 1, \\
    c_0 A_0 + c_1 A_1 = 0,
\end{cases} \]
with the constants $A_k$ for $k \in \{0, 1\}$ defined by
\[ A_k = \int_0^{2\:\pi} \psi_{1, k} \psi_0 \, dx. \]
Here, the first condition means that
\[ \int_0^{2\:\pi} f^2(x) \, dx = c_0^2 \int_0^{2\:\pi} \psi_{1, 0}^2(x) \, dx + c_1^2 \int_0^{2\:\pi} \psi_{1, 1}^2(x) \, dx = c_0^2 + c_1^2 = 1 \]
while the second condition implies
\[ f_0 = \int_0^{2\:\pi} f(x) \psi_0(x) \, dx = c_0 A_0 + c_1 A_1 = 0.\]
Applying \eqref{eq:dirichlet_integral_eigenfunc_orthonormal} to $J_1$, we have
\[ J_1(f, f) = c_0^2 \lambda_{1, 0} + c_1^2 \lambda_{1, 1} \leq \lambda_{1, 1} \]
while applying Lemma \ref{lem:dirichlet_integral_parseval} along with the fact that $f_0 = 0$ gives
\[ J(f, f) \geq \sum_{n=1}^\infty \lambda_n f^2_n \geq \lambda_1 \sum_{n=1}^\infty f^2_n = \lambda_1 \int_0^{2\:\pi} f^2(x) \, dx = \lambda_1. \]
Combining these two statements gives
\[ \lambda_{1, 1} \geq J_1(f, f) \geq J(f, f) \geq \lambda_1 \]
which completes the case for $n = 1$. The same approach applies for a general natural number $n$ - we consider
\[ f(x) = \sum_{j=0}^n c_j \psi_{1, j} \]
along with the $n$ homogeneous linear equations given by
\[ f_j = 0 \]
for all $0 \leq j \leq n-1$. Following the same method as in the proof for $n = 1$ shall be sufficient to finish the result. 
\end{proof}

\subsection{Roots of eigenfunctions}

In this section, we aim to identify the roots of the eigenfunctions of \eqref{eq:eigenvalue_relation}, subject to the boundary conditions \eqref{eq:eigenfunc_periodic_boundary_condition}. To do this, we first review the well-known Sturm Comparison Theorem, which illustrates the oscillatory nature of the solutions as well as allowing us to estimate the number of roots. We shall use the following formulation of Simons \cite{simon_2005}.

\begin{thm}[Sturm Comparison Theorem] \label{thm:sturm_comparison}
    Let $V$ be a bounded function and let $\lambda_1$ and $\lambda_2$ be real numbers. For $i \in \{1, 2\}$, let $\psi_i$ be a non-trivial solution to the corresponding version of \eqref{eq:eigenvalue_relation},
    \[ \lambda_i \psi = - \psi'' + V \psi. \]
    Furthermore, suppose that $\psi_1$ has two roots $a < b$. Then $\psi_2$ has a root in $(a, b)$ if one of the following two conditions hold:
    \begin{itemize}
        \item $\lambda_1 < \lambda_2$, or
        \item $\lambda_1 = \lambda_2$ and $\psi_2(a) \neq 0$.
    \end{itemize}
\end{thm}

We start by denoting by $\mu_n$ and $\xi_n$ the eigenvalues and eigenfunctions of \eqref{eq:eigenvalue_relation} with the semi-periodic boundary condition
\begin{align} \label{eq:semi_periodic_boundary_condition}
    \begin{cases}
        \xi(0) = -\,\xi(2\:\pi), \\
        \xi'(0) = -\,\xi'(2\:\pi).
    \end{cases}
\end{align}
and $\Lambda_n$ and $\Psi_n$ the eigenvalues and eigenfunctions of \eqref{eq:eigenvalue_relation} with the Dirichlet boundary condition
\begin{equation} \label{eq:homogeneous_boundary_condition}
    \Psi(0) = \Psi(2\:\pi) = 0.
\end{equation}

Here, the number of roots of the eigenfunctions $\Psi_n$ is well-established in the literature on Sturm-Liouville operators; see the following theorem in Eastham \cite{eastham_ode_1970}.
\begin{thm} \label{thm:zeros_eigenfunction_homogeneous_boundary}
     The eigenfunction $\Psi_n$ has exactly $n$ roots in the open interval $(0, 2\:\pi)$.
\end{thm}

We shall also utilize the following theorem in Eastham \cite{eastham_1973} to compare the eigenvalues of the associated boundary value problems. 

\begin{thm} \label{thm:interlacing_periodic_vs_semiperiodic}
    The eigenvalues $\lambda_n$ are interlaced with the eigenvalues $\mu_n$ of \eqref{eq:eigenvalue_relation} with boundary conditions \eqref{eq:semi_periodic_boundary_condition} according to
    \begin{equation} \label{eq:interlacing_periodic_vs_semiperiodic}
        \lambda_0 < \mu_0 \leq \mu_1 < \lambda_1 \leq \lambda_2 < \mu_2 \leq \mu_3 < \lambda_3 \leq \lambda_4 < ...
    \end{equation}
\end{thm}

\begin{thm} \label{thm:interlacing_periodic_semiperiodic_vs_dirichlet}
    The eigenvalues $\lambda_n$ of \eqref{eq:eigenvalue_relation} with boundary conditions \eqref{eq:eigenfunc_periodic_boundary_condition} and the eigenvalues $\Lambda_n$ of \eqref{eq:eigenvalue_relation} with boundary conditions \eqref{eq:homogeneous_boundary_condition} are related by
    \begin{equation} \label{eq:eigenvalue_periodic_vs_dirichlet}
        \lambda_{2m+1} \leq \Lambda_{2m+1} \leq \lambda_{2m+2}
    \end{equation}
    for any $m \geq 0$. On the other hand, the eigenvalues $\Lambda_n$ are related to the eigenvalues $\mu_n$ of \eqref{eq:eigenvalue_relation} with boundary conditions \eqref{eq:semi_periodic_boundary_condition} by
    \begin{equation} \label{eq:eigenvalue_semiperiodic_vs_dirichlet}
        \mu_{2m} \leq \Lambda_{2m} \leq \mu_{2m+1}.
    \end{equation}
\end{thm}

An immediate corollary of \eqref{eq:interlacing_periodic_vs_semiperiodic}, \eqref{eq:eigenvalue_periodic_vs_dirichlet}, and \eqref{eq:eigenvalue_semiperiodic_vs_dirichlet} is the following interlacing of $\lambda_n$ and $\Lambda_n$.
\begin{equation} \label{eq:interlacing_periodic_vs_dirichlet}
    \lambda_0 < \Lambda_0 < \lambda_1 \leq \Lambda_1 \leq \lambda_2 < \Lambda_2 < \lambda_3 \leq \Lambda_3 \leq \lambda_4 < ...
\end{equation}

These results shall be utilized in proving the following key theorem regarding the roots of the eigenfunctions of our original problem.

\begin{thm} \label{thm:zeros_eigenfunction_periodic_boundary}
Let $\psi_n$ be the eigenfunctions of \eqref{eq:eigenvalue_relation} with boundary conditions \eqref{eq:eigenfunc_periodic_boundary_condition}.  Then
\begin{enumerate}[(i)] 
    \item $\psi_0$ has no roots in $[0, 2\:\pi]$;
    \item For all $m \geq 0$, $\psi_{2m+1}$ and $\psi_{2m+2}$ have precisely $2m+2$ roots in $[0, 2\:\pi)$.
\end{enumerate}
\end{thm}

\begin{proof}
For part {\it(i)}, by Theorem \ref{thm:zeros_eigenfunction_homogeneous_boundary}, the eigenfunction $\Psi_0$ has no roots in $(0, 2\:\pi)$. Furthermore, by \eqref{eq:interlacing_periodic_vs_dirichlet}, $\lambda_0 < \Lambda_0$, so by the Sturm Comparison Theorem, $\psi_0$ cannot have two roots in $[0, 2\:\pi]$. However, because of its periodic boundary conditions \eqref{eq:eigenfunc_periodic_boundary_condition}, $\psi_0$ must have an even number of roots in $[0, 2\:\pi)$. Hence, $\psi_0$ has no roots in $[0, 2\:\pi]$.

For part {\it(ii)}, consider $\psi_{2m+1}$ for any $m \geq 0$. From \eqref{eq:interlacing_periodic_vs_dirichlet}, we find
\[ \Lambda_{2m} < \lambda_{2m+1} \leq \Lambda_{2m+1}. \]
Here, note that by Theorem \ref{thm:zeros_eigenfunction_homogeneous_boundary}, the eigenfunctions $\Psi_{2m}$ and $\Psi_{2m+1}$ have $2m$ and $2m+1$ roots, respectively. Thus, the Sturm Comparison Theorem suggests that $\psi_{2m+1}$ can have at least $2m+1$ and at most $2m+2$ roots. However, again because of its periodic boundary conditions \eqref{eq:eigenfunc_periodic_boundary_condition}, $\psi_{2m+1}$ needs to have an even number of roots in $[0, 2\:\pi)$, which means that its number of roots in this interval is $2m+2$. The same argument applies to the eigenfunction $\psi_{2m+2}$, but with the relation
\[ \Lambda_{2m+1} \leq \lambda_{2m+2} < \Lambda_{2m+2}. \]
This completes the proof.
\end{proof}

\section{Eigenvalue and eigenfunction asymptotics}

\subsection{The Pr\"ufer transformation}

In his 1926 paper, Pr\"ufer \cite{prufer_1926} introduced a transformation that greatly simplifies  the analysis of Sturm-Liouville operators, as it allows the utilization of polar coordinates to obtain existence of eigenvalues and oscillation of eigenfunctions. Consider a differential equation with positive real-valued coefficients $C(x)$ and $D(x)$ with piecewise continuous derivatives
\begin{equation} \label{eq:prufer_general_form}
    [C(x) u'(x)]' + D(x) u(x) = 0 \qquad \roq{on} x \in [x_1, x_2].
\end{equation}
Set 
\[ R(x) = [C(x) D(x)]^{1/2}. \]
For every non-trivial, real-valued solution $u(x)$ of \eqref{eq:prufer_general_form}, we can transfer to polar coordinates by writing
\begin{equation} \label{eq:prufer_transform}
    R(x) u(x) = \rho(x) \sin \theta(x), \qquad \qquad C(x) u'(x) = \rho(x) \cos \theta(x), 
\end{equation}
so that
\begin{equation} \label{eq:prufer_transform_variable}
    \begin{cases}
        \rho = [R^2 u^2 + C^2 (u')^2]^{1/2}, \\
        \tan \theta = Ru / Cu'.
    \end{cases}
\end{equation}
We complete the definition by selecting a point $a_0 \in [x_1, x_2]$ to specify that
\[ -\pi < \theta(a_0) < \pi \]
and, in particular, if $u(a_0) \geq 0$, then
\begin{equation} \label{eq:positive_phase_condition}
    0 \leq \theta(a_0) < \pi.
\end{equation}
With these definitions, the following two properties --- the formula for the derivative of the phase function $\theta(x)$, and the estimate for the value of the phase function depending on the number of zeros --- will be extremely important for our asymptotics.
\begin{pro} \label{pro:phase_derivative}
    The derivative of the phase function $\theta'(x)$ is given by the formula
    \begin{equation} \label{eq:phase_derivative}
        \theta'(x) = \left(\frac{D(x)}{C(x)} \right)^{1/2} + \frac{(C(x)D(x))'}{4 C(x) D(x)} \sin(2\:\theta(x)).
    \end{equation}
\end{pro}

\begin{proof}
    Differentiating the first equation in \eqref{eq:prufer_transform} gives
    \[ R'(x) u(x) + R(x) u'(x) = \rho'(x) \sin(\theta(x)) + \rho(x) \cos(\theta(x)) \theta'(x), \]
    which, after rearrangement, can be evaluated as
    \begin{align} \label{eq:phase_derivative_proof_1}
        \cos(\theta(x)) \theta'(x) &= -\frac{\rho'(x)}{\rho(x)} \sin(\theta(x)) + \frac{R'(x)}{\rho(x)} u(x) + \frac{R(x)}{\rho(x)} u'(x) \notag \\
        &= -\frac{\rho'(x)}{\rho(x)} \sin(\theta(x)) + \frac{R'(x)}{R(x)} \sin(\theta(x)) + \frac{R(x)}{C(x)} \cos(\theta(x))
    \end{align}
    where on the last line we have substituted $u(x)$ and $u'(x)$ using the appropriate equation from \eqref{eq:prufer_transform}. On the other hand, by substituting $u(x)$ and $u'(x)$ using the same equations into \eqref{eq:prufer_general_form}, we obtain
    \[ 0 = D(x) \frac{\rho(x) \sin(\theta(x))}{R(x)} + \rho'(x) \cos(\theta(x)) - \rho(x) \sin(\theta(x)) \theta'(x) \]
    which translates to
    \begin{equation} \label{eq:phase_derivative_proof_2}
        \sin(\theta(x)) \theta'(x) = \frac{\rho'(x)}{\rho(x)} \cos(\theta(x)) + \frac{D(x)}{R(x)} \sin(\theta(x)).
    \end{equation}
    Finally, from equations \eqref{eq:phase_derivative_proof_1} and \eqref{eq:phase_derivative_proof_2}, we eliminate the term containing $\rho'(x)/\rho(x)$ to obtain
    \begin{align*}
        \theta'(x) &= [\sin^2(\theta(x)) + \cos^2(\theta(x))] \theta'(x) = \frac{R'(x)}{R(x)} \sin(\theta(x)) \cos(\theta(x)) + \frac{R(x)}{C(x)} \cos^2(\theta(x)) + \frac{D(x)}{R(x)} \sin^2(\theta(x)) \\
        &= \frac{(C(x)D(x))'}{4 C(x) D(x)} \sin(2\:\theta(x)) + \left(\frac{D(x)}{C(x)} \right)^{1/2}
    \end{align*}
    where in the last equality we substituted $R(x) = (C(x) D(x))^{1/2}$. This proves the identity \eqref{eq:phase_derivative}.
\end{proof}

\begin{pro} \label{pro:phase_range}
    If a solution $u(x)$ of \eqref{eq:prufer_general_form} has $N$ roots in $(a_0, y]$ with $u(a_0) \geq 0$ and $a_0 < y \leq x_2$, then
    \begin{equation} \label{eq:phase_range}
        N\:\pi \leq \theta(y) < (N+1)\:\pi
    \end{equation}
\end{pro}

\begin{proof}
    First, note that by the first equation of \eqref{eq:prufer_transform}, $u(x)$ is zero only when $\theta(x)$ is a multiple of $\pi$. Applying this to the derivative formula \eqref{eq:phase_derivative} shows that $\theta'(x)$ is positive at any root of $u(x)$. Since $u(x)$ and $u'(x)$ are never simultaneously zero, we know that $\rho(x)$ is always positive. We can then refine the previous statement to say that $u(x)$ is zero if and only if $\theta(x)$ is a multiple of $\pi$.
    
    Let the roots of $u(x)$ in $(a_0, y]$ be
    \[ a_0 < \alpha_1 < ... < \alpha_N \leq y \]
    and consider first the interval $(a_0, \alpha_1]$. From the assumptions, we know that 
    \begin{itemize}
        \item $0 \leq \theta(a_0) < \pi$,
        \item $\theta(x)$ is not a multiple of $\pi$ anywhere in $(a_0, \alpha_1)$ since $u(x)$ has no roots there, and 
        \item $\theta(\alpha_1) > 0$.
    \end{itemize}
    These facts allow us to conclude that $\theta(\alpha_1) = \pi$. We shall proceed accordingly --- for the interval $(\alpha_1, \alpha_2]$, note that
    \begin{itemize}
        \item $\theta(x)$ is not a multiple of $\pi$ anywhere in $(\alpha_1, \alpha_2)$ since $u(x)$ has no roots here, and 
        \item $\theta(\alpha_2) > 0$.
    \end{itemize}
    We conclude from these that $\theta(\alpha_2) = 2\:\pi$. Continuing this process shall give us the desired inequality \eqref{eq:phase_range}, where the inequality on the left occurs if $\alpha_N = y$.
\end{proof}

Our goal is to apply this transformation to the Schr\"odinger eigenvalue equation \eqref{eq:eigenvalue_relation}. Because of the conditions imposed on $C(x)$ and $D(x)$, we consider differential equations of the form
\[ u'' + (\lambda - V_1)u = 0 \]
where $V_1$ is a continuously differentiable function of period $2\:\pi$. Since we are concerned with the estimates as $\lambda \to \infty$, the second coefficient $\lambda - V_1$ can also be considered positive for $\lambda$ sufficiently large. The phase function, under these assumptions, becomes a two-parameter function that, since $u(x)$ has period $2\:\pi$, satisfies
\begin{equation} \label{eq:phase_periodicity}
    \theta(2\:\pi, \lambda) - \theta(0, \lambda) = 2k\pi
\end{equation}
for some integer $k$. Furthermore, following \eqref{eq:phase_derivative}, it has the derivative
\begin{equation} \label{eq:phase_derivative_eigenvalue}
    \theta'(x, \lambda) = (\lambda - V_1(x))^{1/2} - \frac{V'_1(x)}{4(\lambda - V_1(x))} \sin(2\:\theta(x, \lambda)).
\end{equation}

\subsection{Eigenvalue asymptotics}

The proof of Theorem \ref{thm:revival} relies on asymptotic estimates of eigenvalues and eigenfunctions --- in particular, whether the eigenvalues are sufficiently close to the square of an integer and the eigenfunctions are sufficiently close to the trigonometric basis of a Fourier series. Both of these properties have been observed in numerical simulations, even in the presence of discontinuities in the potential function. In particular, we deduce the following estimate for the smaller class of $C^2$ potential functions, in which the asymptotic eigenvalue estimate is similar to that of Eastham \cite{eastham_1973}.

\begin{thm} \label{thm:eigenvalue_asymptotic}
For a fixed $M > 0$, consider the differential equation \eqref{eq:eigenvalue_relation} with a real-valued, $2\pi$-periodic potential $V_1(x)$ with continuous derivative up to the second order satisfying
\begin{equation} \label{eq:potential_inf_bound}
    ||V_1||_{L^\infty[0, 2\pi]} \leq M.
\end{equation}
Let $\lambda_{1, n}$ be the ascending eigenvalues of the problem, along with the corresponding eigenfunctions $\psi_{1, n}(x)$. Then, there exist constants $C_M > 0$ and $m_M \in \mathbb{N}$ depending only on $M$ such that for all $m \geq m_M$, the eigenvalues $\lambda_{1, 2m+j}$ with $j \in \{1, 2\}$ satisfy
\begin{equation} \label{eq:eigenvalue_uniform_estimate}
    \left| \lambda_{1, 2m+j}^{1/2} - \left[ (m+1) + \frac{A_{1, 1}}{m+1} \right] \right| \leq C_M \left[ (m+1)^{-3} \left( 1 + ||V''_1||_{L^1[0, 2\pi]} \right) + (m+1)^{-5} ||V'_1||^2_{L^2[0, 2\pi]} \right],
\end{equation}
where the coefficient 
\begin{align} \label{eq:eigenvalue_asymptotic_coefficients}
    A_{1, 1} = \frac{1}{4\pi} \int_0^{2\:\pi} V_1(x) \, dx
\end{align}
is independent of $m$. In particular, as $m \to \infty$, we have the asymptotic estimate
\begin{equation} \label{eq:eigenvalue_asymptotic}
    \lambda_{1, 2m+j}^{1/2} = (m+1) + \frac{A_{1, 1}}{(m+1)} + \O(m^{-3}).
\end{equation}
\end{thm}

\begin{proof}
    We shall abbreviate $\lambda_{1, 2m+j}$ for $j \in \{1, 2\}$ as $\lambda$ throughout this proof. First, observe that the eigenvalues with indices $(2m+1)$ and $(2m+2)$ of \eqref{eq:eigenvalue_relation} when the potential is constant, specifically $V(x) = \pm M$, are $(m+1)^2 \pm M$. Therefore, applying Theorem \ref{thm:eigenvalue_variation} to
    \[ -M \leq V_1(x) \leq M \]
    gives that
    \[ (m+1)^2 - M \leq \lambda \leq (m+1)^2 + M. \]
    Now, for $m$ sufficiently large, we get
    \[ (m+1)^2 \geq 4M \]
    and therefore
    \[ \frac{3}{4} (m+1)^2 \leq (m+1)^2 - M \leq \lambda \leq (m+1)^2 + M \leq \frac{5}{4}(m+1)^2. \]
    This means that $\lambda^{1/2} \sim (m+1)$ uniformly across all potentials satisfying \eqref{eq:potential_inf_bound}. Furthermore, we also have the bound
    \begin{equation} \label{eq:eigenvalue_asymptotic_proof_1}
        \lambda - V_1(x) \geq (m+1)^2 - 2M \geq \frac{(m+1)^2}{2}.
    \end{equation}
    Now, without loss of generality, assume that $\psi_{1, n}(0) \geq 0$ for all $n \geq 0$. Here, we apply the Pr\"ufer transform \eqref{eq:prufer_transform_variable} to $u(x) = \psi_{1, 2m+j}(x)$, $j \in \{1, 2\}$ with $a_0 = 0$, which satisfies \eqref{eq:positive_phase_condition}. Equations \eqref{eq:positive_phase_condition} and \eqref{eq:phase_periodicity} then give us the range of the value of the phase function as
    \[ 2k \pi \leq \theta(2\:\pi, \lambda) < (2k+1)\pi \]
    for some integer $k$. However, since Theorem \ref{thm:zeros_eigenfunction_periodic_boundary} dictates that both $\psi_{1, 2m+1}$ and $\psi_{1, 2m+2}$ have $2m+2$ zeros in $(0, 2\:\pi]$, then Property \ref{pro:phase_range} allows us to deduce that
    \[ k = m+1. \]
    We then transform \eqref{eq:phase_periodicity} into
    \begin{align} \label{eq:eigenvalue_asymptotic_proof_2}
        2(m+1)\pi = \theta(2\:\pi, \lambda) - \theta(0, \lambda) \notag &= \int_0^{2\:\pi} \theta'(x, \lambda) \, dx \notag \\
        &= \int_0^{2\:\pi} (\lambda - V_1(x))^{1/2} \, dx - \frac{1}{4} \int_0^{2\:\pi} \frac{V'_1(x)}{\lambda - V_1(x)} \sin(2\:\theta(x, \lambda)) \, dx,
    \end{align}
    where the last equality is derived using \eqref{eq:phase_derivative_eigenvalue}. We shall focus on the latter term first; further substituting using \eqref{eq:phase_derivative_eigenvalue} yields
    \begin{align} \label{eq:eigenvalue_asymptotic_proof_3}
        \frac{1}{4} &\int_0^{2\:\pi} \frac{V'_1(x)}{(\lambda - V_1(x))^{3/2}} \left[ \theta'(x, \lambda) + \frac{V'_1(x)}{4(\lambda - V_1(x))} \sin(2\:\theta(x, \lambda)) \right] \sin(2\:\theta(x, \lambda)) \, dx \notag \\
        &= \frac{1}{8} \int_0^{2\:\pi} \left[ \frac{d}{dx} \frac{V'_1(x)}{(\lambda - V_1(x))^{3/2}} \right] \cos(2\:\theta(x, \lambda)) \, dx + \frac{1}{16} \int_0^{2\:\pi} \frac{(V'_1(x))^2}{(\lambda - V_1(x))^{5/2}} \sin^2(2\:\theta(x, \lambda)) \, dx \notag \\
        &= \frac{1}{8} \int_0^{2\:\pi} \left[ \frac{V''_1(x)}{(\lambda - V_1(x))^{3/2}} + \frac{3}{2} \frac{(V'_1(x))^2}{(\lambda - V_1(x))^{5/2}} \right] \cos(2\:\theta(x, \lambda)) \, dx + \frac{1}{16} \int_0^{2\:\pi} \frac{(V'_1(x))^2}{(\lambda - V_1(x))^{5/2}} \sin^2(2\:\theta(x, \lambda)) \, dx
    \end{align}
    where we used an integration by parts on the term involving $\theta'(x, \lambda)$, and the boundary term vanishes due to the $2\pi$-periodicity of $V_1(x)$. Here, the fact that $|\cos(2 \: \theta(x, \lambda))| \leq 1$ and $|\sin(2 \: \theta(x, \lambda))| \leq 1$, along with \eqref{eq:eigenvalue_asymptotic_proof_1}, gives us that
    \begin{align} \label{eq:eigenvalue_asymptotic_proof_4}
        \frac{1}{4} \int_0^{2\:\pi} \frac{V'_1(x)}{\lambda - V_1(x)} \sin(2\:\theta(x, \lambda)) \, dx \leq C \left(  \frac{||V''_1||_{L^1[0, 2\pi]}}{(m+1)^3} + \frac{||V'_1||^2_{L^2[0, 2\pi]}}{(m+1)^5} \right)
    \end{align}
    where $C$ is a constant depending only on $M$. On the other hand, we have, by the binomial expansion,
    \[ (\lambda - V_1(x))^{1/2} = \lambda^{1/2} - \frac{\lambda^{-1/2}}{2} V_1(x) + R(x, \lambda) \]
    with
    \[ |R(x, \lambda)| \leq C \lambda^{-3/2} |V_1(x)|^2 \]
    for some constant $C$. Therefore, the first term of \eqref{eq:eigenvalue_asymptotic_proof_2} can be expanded as
    \[ \int_0^{2\:\pi} (\lambda - V_1(x))^{1/2} \, dx = 2\:\pi \lambda^{1/2} - 2\pi \lambda^{-1/2} A_{1, 1} + E(\lambda) \]
    where we can further use \eqref{eq:potential_inf_bound} and the fact that $\lambda^{1/2} \sim (m+1)$ to get that the error term $E(\lambda)$ satisfies the bound
    \begin{equation} \label{eq:eigenvalue_asymptotic_proof_5}
        |E(\lambda)| \leq C \lambda^{-3/2} ||V_1||^2_{L^2[0, 2\pi]} \leq (C \cdot 2\pi M^2) (m+1)^{-3}.
    \end{equation}
    With \eqref{eq:eigenvalue_asymptotic_proof_4} and \eqref{eq:eigenvalue_asymptotic_proof_5}, we get that for some constant $C$ dependent only on $M$,
    \begin{align*}
        \Bigg \vert 2(m+1)\pi - \left( 2\:\pi \lambda^{1/2} - 2\pi \lambda^{-1/2} A_{1, 1} \right) \Bigg \vert \leq C \left[ (m+1)^{-3} \left( 1 + ||V''_1||_{L^1[0, 2\pi]} \right) + (m+1)^{-5} ||V'_1||^2_{L^2[0, 2\pi]} \right].
    \end{align*}
    Here, by absorbing the $2\pi$ factor into the constant $C$, we get
    \begin{align} \label{eq:eigenvalue_asymptotic_proof_6}
        \left| \lambda^{1/2} - \left[ (m+1) + \lambda^{-1/2} A_{1, 1} \right] \right| \leq C \left[ (m+1)^{-3} \left( 1 + ||V''_1||_{L^1[0, 2\pi]} \right) + (m+1)^{-5} ||V'_1||^2_{L^2[0, 2\pi]} \right].
    \end{align}
    Finally, note that $|A_{1, 1}| \leq M/2$, so
    \begin{align} \label{eq:eigenvalue_asymptotic_proof_7}
        \left| \frac{A_{1, 1}}{m+1} - \lambda^{-1/2} A_{1, 1} \right| = |A_{1, 1}| \left| \frac{\lambda^{1/2} - (m+1)}{\lambda^{1/2} (m+1)} \right| &\leq \frac{M}{2} \left| \frac{|\lambda - (m+1)^2|}{\lambda^{1/2} (m+1) [\lambda^{1/2} + (m+1)]} \right| \notag \\
        &\leq \frac{M^2}{2 \lambda^{1/2} (m+1)^2} \notag \\
        &\leq \frac{C}{(m+1)^3}
    \end{align}
    where the constant $C$ only depends on $M$. Combining \eqref{eq:eigenvalue_asymptotic_proof_6} and \eqref{eq:eigenvalue_asymptotic_proof_7} then gives us the required estimate \eqref{eq:eigenvalue_uniform_estimate} for a merged constant $C_M$. The eigenvalue asymptotic \eqref{eq:eigenvalue_asymptotic} follows naturally by sending $m \to \infty$.
\end{proof}

Here, notice that our potential given by Assumption \ref{ass:potential} is only piecewise continuously differentiable up to the second order, so \eqref{eq:eigenvalue_asymptotic} does not apply directly. However, thanks to Theorem \ref{thm:eigenvalue_variation}, we can still formulate a weaker version of this asymptotic result. The approach shall be to bracket the piecewise $C^2$ potential function between two $C^2$ potential functions, on which we can apply Theorem \ref{thm:eigenvalue_asymptotic} to get eigenvalue asymptotics, and then derive the remaining error from there. To achieve this, we first introduce some estimates for the bracketing functions in the following lemma.

\begin{lem} \label{lem:bracket_estimate}
    Consider a real-valued piecewise $C^2$ potential $V(x)$ as defined in Assumption \ref{ass:potential}. Then, there exists $\delta_0 \in (0, 1)$ such that for any \textit{smoothing width} $0 < \delta \leq \delta_0$, there exist real-valued, $2\pi$-periodic $C^2$ functions $V^+(x)$ and $V^-(x)$ such that
    \begin{equation} \label{eq:bracket_estimate_monotonicity}
        V^+(x) \geq V(x) \geq V^-(x) \qquad \roq{for a.e.} x \in [0, 2\pi]
    \end{equation}
    and there exists $M > 0$ such that
    \begin{equation} \label{eq:bracket_estimate_inf_bound}
        \max\{ ||V^+||_{L^\infty[0, 2\pi]}, ||V^-||_{L^\infty[0, 2\pi]} \} \leq M.
    \end{equation}
    Moreover, the functions can be chosen so that $V^{\pm}(x)$ is close to $V(x)$ in the $L^1$- and $L^2$-norms, given by
    \begin{equation} \label{eq:bracket_estimate_norm_closeness}
        ||V^{\pm} - V||_{L^1[0, 2\pi]} \leq C \delta \quad \roq{and} \quad ||V^{\pm} - V||_{L^2[0, 2\pi]} \leq C \delta^{1/2}.
    \end{equation}
    Furthermore, we have, on the first derivative of $V^{\pm}(x)$, the estimate
    \begin{equation} \label{eq:bracket_estimate_first_derivative}
        ||(V^{\pm})'||^2_{L^2[0, 2\pi]} \leq C \left( 1+ \delta^{-1} \right),
    \end{equation}
    and, on the second derivative of $V^{\pm}(x)$, the estimate
    \begin{equation} \label{eq:bracket_estimate_second_derivative}
        ||(V^\pm)''||_{L^1[0, 2\pi]} \leq C \left( 1+ \delta^{-1} \right).
    \end{equation}
    In addition, with $\{x_k\}$ being the discontinuities of $V$ as in Assumption \ref{ass:potential}, if we define, for $0 \leq n \leq 2$ and $0 \leq k \leq N$,
    \begin{equation} \label{eq:discontinuity_jumps}
        h_{k, n} := \lim_{x \to x_k^+} V^{(n)}(x) - \lim_{x \to x_k^-} V^{(n)}(x),
    \end{equation}
    then the constants above depend only on $N$, the values $\{h_{k, n}\}$, and the norms $||V'||_{L^2[0, 2\pi]}$,  $||V''||_{L^1[0, 2\pi]}$, and $||V||_{L^\infty[0, 2\pi]}$. In particular, they are uniform across all $\delta \in (0, \delta_0]$.
    
\end{lem}

\begin{proof}
    Within this proof, the discontinuities are interpreted periodically on $\RR$ with period $2\pi$. First, we choose
    \[ \delta_0 = \min \left\{ \frac{1}{2}, \frac{1}{4} \min_{0 \leq k \leq N} \{|x_{k+1} - x_{k}|\} \right\} , \]
    which guarantees that the smoothing width \(\delta\) has to be smaller than one quarter of the minimum distance between consecutive discontinuity points, so the smoothing neighborhoods are pairwise disjoint. Then, we form a function $V^+(x)$ by adding, at the lower side of each discontinuity, a positive smoothing function over a region of width $\delta$ so that the resulting function is $C^2$. Specifically, let $H_0(x)$, $H_1(x)$, $H_2(x)$ be polynomials such that
    \[ H^{(n)}_0(0) = H^{(n)}_1(0) = H^{(n)}_2(0) = 0 \]
    for all $0 \leq n \leq 2$, while
    \begin{align*}
        H_0(1) = 1, \quad &H_1(1) = 0, \quad H_2(1) = 0, \\
        H'_0(1) = 0, \quad &H'_1(1) = 1, \quad H'_2(1) = 0, \\
        H''_0(1) = 0, \quad &H''_1(1) = 0, \quad H''_2(1) = 1.
    \end{align*}
    Explicitly, we can take 
    \begin{align*}
        H_0(x) &= 6x^5 - 15x^4 + 10x^3, \\
        H_1(x) &= -3x^5 + 7x^4 - 4x^3, \\
        H_2(x) &= \frac{1}{2} x^5 - x^4 + \frac{1}{2} x^3.
    \end{align*}
    Then, using the sign function convention
    \[ \sgn(x) = \begin{cases}
        1 \quad &\roq{if} x \geq 0, \\
        -1 &\roq{if} x < 0,
    \end{cases} \]
    and with the values $\{h_{k, n}\}$ as specified in \eqref{eq:discontinuity_jumps}, we define the function $\eta_k$ by
    \begin{equation*}
        \eta_k(x) = \begin{cases}
            |h_{k, 0}| H_0(y) + \delta h_{k, 1} H_1(y) + \sgn(h_{k, 0}) \cdot \delta^2 h_{k, 2} H_2(y) \quad &\roq{for} x \in I_{k, \delta}, \\
            0 &\roq{otherwise,}
        \end{cases}
    \end{equation*}
    where
    \[ y = 1 + \frac{\sgn(h_{k, 0}) \cdot (x - x_k)}{\delta} \]
    and $I_{k, \delta}$, being the smoothing region, is defined as
    \begin{equation*}
        I_{k, \delta} = \begin{cases}
            [x_k - \delta, x_k) \quad &\roq{if} h_{k, 0} \geq 0, \\
            (x_k, x_k + \delta] &\roq{if} h_{k, 0} < 0.
        \end{cases}
    \end{equation*}
    Furthermore, we define a $C^2$ bump
    \[ B(x) = \begin{cases}
        (1 - x^2)^3 \quad &\roq{for all} x \in [-1, 1], \\
        0 &\roq{otherwise,}
    \end{cases} \]
    which satisfies $B^{(n)}(\pm 1) = 0$ for all $0 \leq n \leq 2$, and let
    \[ C_k = \max \left\{ 4\delta|h_{k, 1}| + \frac{\delta^2}{2} |h_{k, 2}| - |h_{k, 0}|, \, 0 \right\}. \]
    Upon checking, one can verify that the function $V + \eta_k$ has continuous derivatives up to the second order at $x_k$, and hence
    \begin{equation} \label{eq:bracket_estimate_proof_1}
        V^+(x) = V(x) + \sum_{k = 0}^N \left[ \eta_k(x) + C_k B\left( \frac{x - x_k}{\delta} \right) \right]
    \end{equation}
    is a $C^2$ function. Note that the formula above only defines the value of $V^+$ away from the discontinuity points, so we define the value of $V^+$ at $x_k$ by continuous extension, which is possible because the one-sided values and first two derivatives match. Furthermore, at $x_0 = 0$, we interpret the function periodically - that is, its left side is $x \to 2\pi^-$ and its right side is $x \to 0^+$. 
    
    We now start verifying the required properties. First, note that $x \in I_{k, \delta}$ corresponds to $y \in [0, 1)$ under the change of variables. Now, we have the bounds
    \begin{align*}
        y^{-3} H_0(y) = 6y^2 - 15y + 10 \geq 1 \qquad \roq{for all} y \in [0, 1)
    \end{align*}
    while, for $y \in [0, 1)$, 
    \begin{align*}
        |y^{-3} H_1(y)| &= |3y^2 - 7y + 4| \leq 4, \\
        |y^{-3} H_2(y)| &= \left| \frac{1}{2} y^2 - y + \frac{1}{2} \right| \leq \frac{1}{2},
    \end{align*}
    This means that, for $x \in I_{k, \delta}$,
    \begin{align*}
        \eta_k(x) + C_k B\left( \frac{x - x_k}{\delta} \right) &\geq y^3 \left( |h_{k, 0}| - 4\delta|h_{k, 1}| - \frac{\delta^2}{2} |h_{k, 2}| \right) + C_k \cdot y^3(2-y)^3 \\
        &\geq y^3 \left( |h_{k, 0}| - 4\delta|h_{k, 1}| - \frac{\delta^2}{2} |h_{k, 2}| + C_k \right) \\
        &\geq 0.
    \end{align*}
    where it becomes apparent that our choice of $C_k$ above is to guarantee that the right-hand side is non-negative. Then,
    \begin{equation*}
        V^+(x) = V(x) + \sum_{k = 0}^N \left[ \eta_k(x) + C_k B\left( \frac{x - x_k}{\delta} \right) \right] \geq V(x) \qquad \roq{for a.e.} x \in [0, 2\pi],
    \end{equation*}
    which gives us one side of \eqref{eq:bracket_estimate_monotonicity}. Next, we have that, for $y \in [0, 1)$, 
    \begin{equation*}
        |H_0(y)| \leq 10, \quad |H_1(y)| \leq 4, \roq{and} |H_2(y)| \leq \frac{1}{2},
    \end{equation*}
    so we have
    \begin{align*}
        ||V^+ - V||_{L^1[0, 2\pi]} &\leq \sum_{k = 0}^N \left[ \int_{x_k - \delta}^{x_k +\delta} \left| C_k B\left( \frac{x - x_k}{\delta} \right) \right| dx + \int_{I_{k, \delta}} |\eta_k(x)| \, dx \right] \\
        &\leq \sum_{k = 0}^N \left[ 2C_k \delta + \int_{I_{k, \delta}} 10 |h_{k, 0}| + 4 \delta |h_{k, 1}| + \frac{1}{2} \delta^2 |h_{k, 2}| \, dx \right] \\
        &\leq \delta \sum_{k=0}^N \left( 2C_k + 10 |h_{k, 0}| + 4 \delta |h_{k, 1}| + \frac{1}{2} \delta^2 |h_{k, 2}| \right).
    \end{align*}
    Furthermore, we also have
    \begin{align*}
        ||V^+ - V||_{L^2[0, 2\pi]} &= \left( \sum_{k = 0}^N \int_{x_k - \delta}^{x_k + \delta} \left( \eta_k(x) + C_k B\left( \frac{x - x_k}{\delta} \right) \right)^2 dx \right)^{1/2} \\
        &\leq \left[ 2 \sum_{k = 0}^N \left( \int_{I_{k, \delta}} \eta^2_k(x) \, dx + \int_{x_k - \delta}^{x_k + \delta} \left( C_k B\left( \frac{x - x_k}{\delta} \right) \right)^2 dx \right) \right]^{1/2} \\
        &\leq \left[ 2 \sum_{k = 0}^N \left( \int_{I_{k, \delta}} 3 \left( 100 h^2_{k, 0} + 16 \delta^2 h^2_{k, 1} + \frac{1}{4} \delta^4 h^2_{k, 2} \right) dx + 2 \delta C_k^2 \right) \right]^{1/2} \\
        &\leq \left[ 2 \delta \sum_{k=0}^N \left( 300 h^2_{k, 0} + 48 \delta^2 h^2_{k, 1} + \frac{3}{4} \delta^4 h^2_{k, 2} + 2 C_k^2 \right) \right]^{1/2} \\
        &= \sqrt{2\delta} \left[ \sum_{k=0}^N \left( 300 h^2_{k, 0} + 48 \delta^2 h^2_{k, 1} + \frac{3}{4} \delta^4 h^2_{k, 2} + 2 C_k^2 \right) \right]^{1/2}.
    \end{align*}
    This gives us both estimates in \eqref{eq:bracket_estimate_norm_closeness} where, noticeably, the constants depend only on $N$ and $\{h_{k, n}\}$. Now, differentiating both sides of \eqref{eq:bracket_estimate_proof_1} gives us
    \begin{align*}
        ||(V^+)'||^2_{L^2[0, 2\pi]} &= \int_0^{2\pi} \left( V'(x) + \sum_{k = 0}^N \left( \eta'_k(x) + \frac{C_k}{\delta} B'\left( \frac{x - x_k}{\delta} \right) \right) \right)^2 \, dx \\
        &\leq \int_0^{2\pi} (2N+3) \left( (V'(x))^2 + \sum_{k = 0}^N (\eta'_k(x))^2 + \sum_{k = 0}^N \left( \frac{C_k}{\delta} B'\left( \frac{x - x_k}{\delta} \right) \right)^2 \right) dx \\
        &= (2N+3) \left( ||V'||^2_{L^2[0, 2\pi]} + \sum_{k = 0}^N \int_{I_{k, \delta}} (\eta'_k(x))^2 \, dx + \sum_{k = 0}^N \frac{C_k^2}{\delta^2} \int_{x_k - \delta}^{x_k + \delta} \left( B'\left( \frac{x - x_k}{\delta} \right) \right)^2 \, dx \right),
    \end{align*}
    where since we have that, for $y \in [0, 1)$, 
    \begin{equation*}
        |H'_0(y)| \leq 30, \quad |H'_1(y)| \leq 12, \roq{and} |H'_2(y)| \leq \frac{3}{2},
    \end{equation*}
    then
    \begin{align*}
        \int_{I_{k, \delta}} (\eta'_k(x))^2 \, dx &\leq \int_{I_{k, \delta}} \left( \frac{1}{\delta}|h_{k, 0} H'_0(y)| + |h_{k, 1} H'_1(y)| + \delta |h_{k, 2} H'_2(y)| \right)^2 \, dx \\
        &\leq \int_{I_{k, \delta}} \frac{3}{\delta^2} \left( 900|h_{k, 0}|^2 + 144 \delta^2 |h_{k, 1}|^2 + \frac{9\delta^4}{4} |h_{k, 2}|^2 \right) dx \\
        &\leq \frac{3}{\delta} \left( 900|h_{k, 0}|^2 + 144 \delta^2 |h_{k, 1}|^2 + \frac{9\delta^4}{4} |h_{k, 2}|^2 \right)
    \end{align*}
    for all $0 \leq k \leq N$. Therefore, 
    \begin{align*}
        &||(V^+)'||^2_{L^2[0, 2\pi]} \\
        &\leq (2N+3) \left( ||V'||^2_{L^2[0, 2\pi]} + \frac{3}{\delta} \sum_{k=0}^N \left( 900|h_{k, 0}|^2 + 144 \delta^2 |h_{k, 1}|^2 + \frac{9\delta^4}{4} |h_{k, 2}|^2 \right) + \frac{2 ||B'||^2_{L^\infty}}{\delta} \sum_{k=0}^N C_k^2 \right) \\
        &\leq C(1 + \delta^{-1}),
    \end{align*}
    where the last inequality is because $||V'||^2_{L^2[0, 2\pi]}$ is bounded, so we can absorb all constants along with the norm to get \eqref{eq:bracket_estimate_first_derivative}, where $C$ depends only on $||V'||_{L^2}$, $N$, and $\{h_{k, n}\}$. A similar analysis can be performed on the second derivative of $V^+$, as we have
    \begin{align*}
        ||(V^+)''||_{L^1[0, 2\pi]} &= \int_0^{2\pi} \left| V''(x) + \sum_{k = 0}^N \eta''_k(x) \, dx + \sum_{k = 0}^N \frac{C_k}{\delta^2} B''\left( \frac{x - x_k}{\delta} \right) \right| \, dx \\
        &\leq ||V''||_{L^1[0, 2\pi]} + \sum_{k = 0}^N \int_{I_{k, \delta}} \left| \eta''_k(x) \right| dx + \sum_{k = 0}^N \frac{C_k}{\delta^2} \int_{x_k - \delta}^{x_k + \delta} \left| B''\left( \frac{x - x_k}{\delta} \right) \right| dx
    \end{align*}
    where since we have that, for $y \in [0, 1)$, 
    \begin{equation*}
        |H''_0(y)| \leq 60, \quad |H''_1(y)| \leq 24, \roq{and} |H''_2(y)| \leq 3,
    \end{equation*}
    then
    \begin{align*}
        \int_{I_{k, \delta}} \left| \eta''_k(x) \right| dx &\leq \int_{I_{k, \delta}} \frac{1}{\delta^2}|h_{k, 0} H''_0(y)| + \frac{1}{\delta} |h_{k, 1} H''_1(y)| + |h_{k, 2} H''_2(y)| \, dx \\
        &\leq \int_{I_{k, \delta}} \frac{1}{\delta^2} \left( 60|h_{k, 0}| + 24\delta |h_{k, 1}| + 3 \delta^2 |h_{k, 2}| \right) dx \\
        &\leq \frac{1}{\delta} \left( 60|h_{k, 0}| + 24\delta |h_{k, 1}| + 3 \delta^2 |h_{k, 2}| \right)
    \end{align*}
    for all $0 \leq k \leq N$. Therefore,
    \begin{align*}
        ||(V^+)''||_{L^1[0, 2\pi]} &\leq ||V''||_{L^1[0, 2\pi]} + \frac{1}{\delta} \sum_{k=0}^N \left( 60|h_{k, 0}| + 24\delta |h_{k, 1}| + 3 \delta^2 |h_{k, 2}| \right) + \frac{2 ||B''||_{L^\infty}}{\delta} \sum_{k=0}^N C_k \\
        &\leq C(1 + \delta^{-1}),
    \end{align*}
    where the last inequality is because $||V''||_{L^1[0, 2\pi]}$ is bounded, and therefore we can absorb all constants along with the norm to get \eqref{eq:bracket_estimate_second_derivative} with the constant depending only on $||V''||_{L^1}$, $N$, and $\{h_{k, n}\}$. This completes all required estimates for $V^+(x)$. The other bracketing function, $V^-(x)$, will be formed in a similar manner, but now by subtracting a positive smoothing function over a region of width $\delta$ on the higher side of the discontinuity. Specifically, for $0 \leq k \leq N$, with the values $\{h_{k, n}\}$ as specified in \eqref{eq:discontinuity_jumps}, define a function $\zeta_k$ by
    \begin{equation*}
        \zeta_k(x) = \begin{cases}
            |h_{k, 0}| H_0(z) - \delta h_{k, 1} H_1(z) + \sgn(h_{k, 0}) \cdot \delta^2 h_{k, 2} H_2(z) &\roq{for} x \in J_{k, \delta} \\
            0 &\roq{otherwise}
        \end{cases}
    \end{equation*}
    where
    \[ z = 1 - \frac{\sgn(h_{k, 0}) \cdot (x - x_k)}{\delta} \]
    and $J_{k, \delta}$ is defined as
    \begin{equation*}
        J_{k, \delta} = \begin{cases}
            (x_k, x_k + \delta] \quad &\roq{if} h_{k, 0} \geq 0, \\
            [x_k - \delta, x_k) &\roq{if} h_{k, 0} < 0.
        \end{cases}
    \end{equation*}
    Using these smoothing functions, the resulting bracketing function is
    \begin{equation*} 
        V^-(x) = V(x) - \sum_{k = 0}^N \left[ \zeta_k(x) + C_k B\left( \frac{x - x_k}{\delta} \right) \right],
    \end{equation*}
    where the values of $V^-$ at discontinuity points $x_k$ are once again continuously extended from matching both sides, and the value at $x_0 = 0$ is interpreted periodically. By repeating similar techniques as for $V^+(x)$, we can also obtain estimates \eqref{eq:bracket_estimate_monotonicity}, \eqref{eq:bracket_estimate_norm_closeness}, \eqref{eq:bracket_estimate_first_derivative}, and \eqref{eq:bracket_estimate_second_derivative} for $V^-(x)$. It remains for us to prove \eqref{eq:bracket_estimate_inf_bound}. Indeed, 
    \begin{align*}
        ||V^+||_{L^\infty[0, 2\pi]} &\leq ||V||_{L^\infty[0, 2\pi]} + \sum_{k = 0}^N \left( ||\eta_k||_{L^\infty[0, 2\pi]} + C_k ||B||_{L^\infty[0, 2\pi]} \right) \\
        &\leq ||V||_{L^\infty[0, 2\pi]} + \sum_{k = 0}^N \left( 10 |h_{k, 0}| + 4 \delta |h_{k, 1}| + \frac{1}{2} \delta^2 |h_{k, 2}| + C_k \right) \leq M
    \end{align*}
    for some $M$ dependent only on $||V||_{L^\infty}$, $N$, and $\{h_{k, n}\}$. A similar estimate can be applied to $V^-$, after which we can merge the two estimates to get a common constant $M > 0$ satisfying \eqref{eq:bracket_estimate_inf_bound}. This completes the construction.
\end{proof}

The asymptotics for eigenvalues corresponding to piecewise $C^2$ potential functions are now presented in the corollary below.

\begin{thm} \label{cor:eigenvalue_asymptotic_piecewise}
    Consider the differential equation \eqref{eq:eigenvalue_relation} but with the real-valued potential $V(x)$ being piecewise continuously differentiable up to the second order, as defined by Assumption \ref{ass:potential}. Let $\lambda_{n}$ be the ascending eigenvalues of the problem, along with the corresponding eigenfunctions $\psi_{n}(x)$. Then, as $m \to \infty$, the eigenvalues $\lambda_{2m+1}$ and $\lambda_{2m+2}$ satisfy the asymptotic estimate
    \begin{equation} \label{eq:eigenvalue_asymptotic_piecewise}
        \lambda^{1/2} = (m+1) + \frac{A_1}{(m+1)} + \O(m^{-2})
    \end{equation}
    where the coefficient 
    \begin{align} \label{eq:eigenvalue_asymptotic_piecewise_coefficients}
        A_1 = \frac{1}{4\pi} \int_0^{2\:\pi} V(x) \, dx
    \end{align}
    is independent of $m$.
\end{thm}

\begin{proof}
    First, we fix a $\delta_0$ as described in Lemma \ref{lem:bracket_estimate} and form two bracketing functions $V^+(x)$ and $V^-(x)$. Furthermore, let $\lambda^+_n$ be the eigenvalues corresponding to the potential $V^+(x)$, and let $A^+_1$ be defined as in \eqref{eq:eigenvalue_asymptotic_piecewise_coefficients} but with $V^+(x)$ in place of $V(x)$. Similarly, define $\lambda^-_n$ and $A^-_1$ for $V^-(x)$. Using Theorem \ref{thm:eigenvalue_variation}, the fact that
    \begin{equation*}
        V^+(x) \geq V(x) \geq V^-(x) \quad \roq{for a.e.} x \in [0, 2\pi]
    \end{equation*}
    means that for all $n \geq 0$,
    \begin{equation*}
        \lambda^+_n \geq \lambda_n \geq \lambda^-_n.
    \end{equation*}
    On the other hand, \eqref{eq:bracket_estimate_monotonicity} and \eqref{eq:bracket_estimate_norm_closeness} give us
    \[ A_1 \leq A_1^+ \leq A_1 + C\delta. \]
    Furthermore, with the constant $M$ provided in \eqref{eq:bracket_estimate_inf_bound}, we apply Theorem \ref{thm:eigenvalue_asymptotic}; the asymptotic expansion \eqref{eq:eigenvalue_uniform_estimate} gives that, for $n \in \{2m+1, 2m+2\}$,
    \begin{align*}
        \lambda_n^{1/2} &\leq (\lambda^+_n)^{1/2} \\
        &\leq (m+1) + \frac{A^+_1}{m+1} + C \left[ (m+1)^{-3} \left(1 + ||(V^+)''||_{L^1[0, 2\pi]} \right) + (m+1)^{-5} ||(V^+)'||^2_{L^2[0, 2\pi]} \right] \\
        &\leq (m+1) + \frac{A_1}{m+1} + C \left[ (m+1)^{-1} \delta + (m+1)^{-3} \left(1 + ||(V^+)''||_{L^1[0, 2\pi]} \right) + (m+1)^{-5} ||(V^+)'||^2_{L^2[0, 2\pi]} \right],
    \end{align*}
    which means by applying \eqref{eq:bracket_estimate_first_derivative} and \eqref{eq:bracket_estimate_second_derivative}, we get that
    \begin{align*}
        \lambda_n^{1/2} \leq (\lambda^+_n)^{1/2} &\leq (m+1) + \frac{A_1}{m+1} + C \left[ (m+1)^{-1} \delta + (m+1)^{-3} \left(1 + \delta^{-1} \right) + (m+1)^{-5} \left(1 + \delta^{-1} \right) \right] \\
        &\leq (m+1) + \frac{A_1}{m+1} + 2C \left[ (m+1)^{-1} \delta + (m+1)^{-3} \left(1 + \delta^{-1} \right) \right]
    \end{align*}
    for some constant $C$ that now depends on $M$, $||V'||_{L^2}$, $||V''||_{L^1}$, $\{h_{k, n}\}$, and $N$. We can repeat the same process for $V^-(x)$ so that we get
    \begin{align*}
        \lambda_n^{1/2} \geq (\lambda^-_n)^{1/2} \geq (m+1) + \frac{A_1}{m+1} - 2C \left[ (m+1)^{-1} \delta + (m+1)^{-3} \left(1 + \delta^{-1} \right) \right].
    \end{align*}
    Combining the two inequalities above then yields the estimate
    \begin{align*}
        \lambda_n^{1/2} = (m+1) + \frac{A_1}{m+1} + O \left( (m+1)^{-1} \delta + (m+1)^{-3} \left(1 + \delta^{-1} \right) \right),
    \end{align*}
    where the constant of the error term is uniform across all $0 < \delta \leq \delta_0$. Then, we choose $\delta = (m+1)^{-1}$, where for sufficiently large $m$ our estimates apply. This gives
    \begin{align*}
        \lambda_n^{1/2} = (m+1) + \frac{A_1}{m+1} + O \left( m^{-2} \right),
    \end{align*}
    which is our desired estimate. This concludes the proof.
\end{proof}

\subsection{Eigenfunction asymptotics}

Our second goal for this section is to obtain asymptotic estimates for eigenfunctions, partially relying on the asymptotic expansion of eigenvalues \eqref{eq:eigenvalue_asymptotic_piecewise}. To do that, we first acknowledge that any solution $\psi$ of \eqref{eq:eigenvalue_relation} can be written as a linear combination
\begin{equation} \label{eq:fundamental_set_linear_combination}
    \psi(x) = c_1 \phi_1(x) + c_2 \phi_2(x)
\end{equation}
where $\{\phi_1, \phi_2\}$ is a fundamental set of solutions for \eqref{eq:eigenvalue_relation}, satisfying
\begin{equation} \label{eq:fundamental_set}
    \begin{cases}
        \phi_1(0, \lambda) = 1 \\
        \phi'_1(0, \lambda) = 0
    \end{cases} \qquad \qquad \roq{and} \qquad \qquad \begin{cases}
        \phi_2(0, \lambda) = 0 \\
        \phi'_2(0, \lambda) = 1
    \end{cases}.
\end{equation}
Using variation of constants, we can obtain solution representations for $\phi_1$ and $\phi_2$ as
\begin{equation} \label{eq:variation_of_constants_1}
    \phi_1(x, \lambda) = \cos(x\sqrt{\lambda}) + \lambda^{-1/2} \int_0^x \sin((x-y)\sqrt{\lambda}) V(y) \phi_1(y, \lambda) \, dy
\end{equation}
and
\begin{equation} \label{eq:variation_of_constants_2}
    \phi_2(x, \lambda) = \lambda^{-1/2} \sin(x\sqrt{\lambda}) + \lambda^{-1/2} \int_0^x \sin((x-y)\sqrt{\lambda}) V(y) \phi_2(y, \lambda) \, dy.
\end{equation}
Both of these formulas can be regarded as recursive definitions of the fundamental solutions, and therefore we can iterate them as many times as we wish to yield terms of higher order in $\lambda$. Let us present the resulting asymptotic estimates for the two fundamental solutions. 

\begin{thm} \label{thm:eigenfunction_asymptotic}
    Assume that $V(x)$ is a real-valued, piecewise continuous potential function with mean zero:
    \[ \langle V \rangle = \frac{1}{2\:\pi} \int_0^{2\:\pi} V(x) \, dx = 0 \]
    Then, the fundamental solutions \eqref{eq:fundamental_set} satisfy the asymptotics
    \[ \phi_1(x, \lambda) = \cos(x\sqrt{\lambda}) + \lambda^{-1/2} \int_0^x \sin((x-y)\sqrt{\lambda}) \cos(y\sqrt{\lambda}) V(y) \, dy + \O(\lambda^{-1}) \]
    and 
    \[ \phi_2(x, \lambda) = \lambda^{-1/2} \sin(x\sqrt{\lambda}) + \lambda^{-1} \int_0^x \sin((x-y)\sqrt{\lambda}) \sin(y\sqrt{\lambda}) V(y) \, dy + \O(\lambda^{-3/2}). \]
\end{thm}

\begin{proof}
    First, we make some remarks about the boundedness of the fundamental solutions. Let $M_1(\lambda)$ be the maximum absolute value of $\phi_1$ over $[0, 2\:\pi]$. Then, \eqref{eq:variation_of_constants_1} gives
    \[ M_1(\lambda) \leq 1 + \lambda^{-1/2} M_1(\lambda) \int_0^{2\:\pi} |V(y)| \, dy \]
    which means
    \[ M_1(\lambda) \leq \left( 1 - \lambda^{-1/2} \int_0^{2\:\pi} |V(y)| \, dy \right)^{-1} \]
    provided that $\lambda$ is large enough for the right-hand side to be positive. Since $V(x)$ is bounded, 
    \[ M_1(\lambda) = \O(1) \roq{as} \lambda \longrightarrow  \infty.\]
A similar evaluation can be conducted for $\phi_2$ to get that its maximum absolute value over $[0, 2\pi]$ is $\O(\lambda^{-1/2})$ as $\lambda \to \infty$. Here, from \eqref{eq:variation_of_constants_1}, we iterate one more time by replacing the right-hand side in place of $\phi_1(x, \lambda)$ in the integral to get
    \begin{align*}
        \phi_1(x, \lambda) = \cos(x\sqrt{\lambda}) &+ \lambda^{-1/2} \int_0^x \sin((x-y)\sqrt{\lambda}) \cos(y\sqrt{\lambda}) V(y) \, dy \\
        &+ \lambda^{-1} \int_0^x \sin((x-y)\sqrt{\lambda}) \, V(y) \left( \int_0^y \sin((y-z)\sqrt{\lambda}) V(z) \phi_1(z, \lambda) \, dz \right) dy
    \end{align*}
    where, with the upper bound as above, the asymptotic estimate
    \begin{equation} \label{eq:fundamental_sol_1_asymptotic_lambda}
        \phi_1(x, \lambda) = \cos(x\sqrt{\lambda}) + \lambda^{-1/2} \int_0^x \sin((x-y)\sqrt{\lambda}) \cos(y\sqrt{\lambda}) V(y) \, dy + \O(\lambda^{-1})
    \end{equation}
   holds. On the other hand, iterating one more time for \eqref{eq:variation_of_constants_2} gives
    \begin{align*}
        \phi_2(x, \lambda) = \lambda^{-1/2} \sin(x\sqrt{\lambda}) &+ \lambda^{-1} \int_0^x \sin((x-y)\sqrt{\lambda}) \sin(y\sqrt{\lambda}) V(y) \, dy \\
        &+ \lambda^{-1} \int_0^x \sin((x-y)\sqrt{\lambda}) \, V(y) \left( \int_0^y \sin((y-z)\sqrt{\lambda}) V(z) \phi_2(z, \lambda) \, dz \right) dy
    \end{align*}
    and, by substituting the upper bound, 
    \begin{equation} \label{eq:fundamental_sol_2_asymptotic_lambda}
        \phi_2(x, \lambda) = \lambda^{-1/2} \sin(x\sqrt{\lambda}) + \lambda^{-1} \int_0^x \sin((x-y)\sqrt{\lambda}) \sin(y\sqrt{\lambda}) V(y) \, dy + \O(\lambda^{-3/2}).
    \end{equation}
    This completes the proof.
\end{proof}

Here, note that with the mean zero condition on the potential, \eqref{eq:eigenvalue_asymptotic_piecewise} implies that the eigenvalues satisfy
\[ \lambda_{2m-1}^{1/2} \simeq \lambda_{2m}^{1/2} = m + \O(m^{-2})\roq{as} m \longrightarrow  \infty. \]
Applying this to \eqref{eq:fundamental_sol_1_asymptotic_lambda} yields
\begin{equation} \label{eq:fundamental_sol_1_asymptotic}
    \phi_1(x, \lambda_{2m-1}) \simeq \phi_1(x, \lambda_{2m}) = \cos(mx) + \frac{1}{m} \int_0^x \sin(m(x-y)) \cos(my) V(y) \, dy + \O(m^{-2}).
\end{equation}
This is because the $\O(m^{-2})$ term within the cosine function translates to $\O(m^{-2})$ maximum error, which is absorbed in the $\O(\lambda^{-1}) = \O(m^{-2})$ term. On the other hand, the $\O(m^{-2})$ term within the sine function also translates to $\O(m^{-2})$ maximum error, so \eqref{eq:fundamental_sol_2_asymptotic_lambda} gives 
\begin{equation} \label{eq:fundamental_sol_2_asymptotic}
    \phi_2(x, \lambda_{2m-1}) \simeq \phi_2(x, \lambda_{2m}) = \frac{1}{m} \sin(mx) + \frac{1}{m^2} \int_0^x \sin(m(x-y)) \sin(my) V(y) \, dy + \O(m^{-3}).
\end{equation}
Finally, knowing that the eigenfunctions $\psi_n$ are linear combinations of the fundamental solutions as in \eqref{eq:fundamental_set_linear_combination}, we have the following lemma on the relative size of the constants.

\begin{lem} \label{lem:coefficient_size}
    For any $m \geq 1$ and $i \in \{2m-1, 2m\}$, the unique coefficients $\alpha_i, \beta_i$ of the linear combination
    \[ \psi_i(x) = \alpha_{i} \phi_1(x, \lambda_i) + \beta_{i} m \phi_2(x, \lambda_i) \]
    satisfy that there exists a constant $C > 0$ independent of $m$ and $i$ such that
    \[ |\alpha_i| + |\beta_i| \leq C. \]
    In particular, both coefficients are $\O(1)$ as $m \to \infty$.
\end{lem}

\begin{proof}
    For $i \in \{2m-1, 2m\}$, \eqref{eq:fundamental_sol_1_asymptotic} and \eqref{eq:fundamental_sol_2_asymptotic} give that for all $x \in [0, 2\pi]$,
    \[ |\phi_1(x, \lambda_i) - \cos(mx)| \leq Cm^{-1} \]
    and 
    \[ |m\phi_2(x, \lambda_i) - \sin(mx)| \leq Cm^{-1} \]
    where $C$ is independent of $m$ and $i$. However, using 
    \[ ||\cos(mx)||^2_{L^2[0, 2\pi]} = ||\sin(mx)||^2_{L^2[0, 2\pi]} = \pi \]
    plus the fact that $\cos(mx)$ and $\sin(mx)$ are orthogonal, we have
    \begin{align*}
        ||\phi_1(\cdot, \lambda_i)||_2^2 = \pi + O(m^{-1}), \\
        ||m\phi_2(\cdot, \lambda_i)||_2^2 = \pi + O(m^{-1}), \\
        \langle \phi_1(x, \lambda_i), m\phi_2(x, \lambda_i) \rangle = O(m^{-1}).
    \end{align*}
    Therefore, for sufficiently large $m$, we have
    \begin{align*}
        1 = ||\psi_i||^2_{L^2[0, 2\pi]} &= ||\alpha_{i} \phi_1(\cdot, \lambda_i) + \beta_{i} m \phi_2(\cdot, \lambda_i)||^2_{L^2[0, 2\pi]} \\
        &= \pi(\alpha_i^2 + \beta_i^2) + O(m^{-1}(\alpha_i^2 + \beta_i^2)) \\
        &\geq \frac{\pi}{2} (\alpha_i^2 + \beta_i^2)
    \end{align*}
    and therefore, 
    \[ |\alpha_i| + |\beta_i| \leq \sqrt{2(\alpha_i^2 + \beta_i^2)} \leq \frac{2}{\sqrt{\pi}}. \]
    Taking into account the finitely many leftover indices, we see that there exists a $C$ independent of $m$ and $i$ such that
    \[ |\alpha_i| + |\beta_i| \leq C. \]
    We conclude that both coefficients have to be $\O(1)$ as $m \to \infty$.
\end{proof}

From here on, we let \[ a_n = \sqrt{\frac{1}{\pi}} \, \cos(nx) \qquad \roq{and} \qquad b_n = \sqrt{\frac{1}{\pi}} \, \sin(nx) \]
denote the standard orthonormal Fourier basis functions.
With the lemma above, and by adjusting the coefficients $\alpha_i$ and $\beta_i$ to the normalized bases, we can express the eigenfunction $\psi_i$ for $i \in \{2m-1, 2m\}$ as
\begin{align} \label{eq:eigenfunction_asymptotic_2_coeff}
    \psi_i(x) &= \alpha_i a_m(x) + \beta_i b_m(x) + \frac{1}{m} \int_0^x \sin(m(x-y)) [\alpha_i a_m(y) + \beta_i b_m(y)] V(y) \, dy + R_i(x)
\end{align}
where $||R_i||_{L^\infty} = \O(m^{-2})$. With this, the next lemma shall give further information about the constants $\alpha_i$ and $\beta_i$.

\begin{lem} \label{lem:coeff_asymptotic}
    For all $m \geq 1$, there exists $y_m$ such that
    \begin{align*}
        \alpha_{2m-1} &= \cos(y_m) + \O(m^{-1}), \hskip40pt \alpha_{2m} = \sin(y_m) + \O(m^{-1}),\\
        \beta_{2m-1} &= \sin(y_m) + \O(m^{-1}), \hskip40pt \beta_{2m} = -\cos(y_m) + \O(m^{-1}).
    \end{align*}
\end{lem}

\begin{proof}
Using \eqref{eq:eigenfunction_asymptotic_2_coeff} and the fact that $\psi_{2m-1}$ is normalized, we have
\[ 1 = ||\psi_{2m-1}||^2_{L^2[0, 2\:\pi]} = \alpha_{2m-1}^2 + \beta_{2m-1}^2 + \O(m^{-1}).\]
This means that there exists a $y_m$ such that
\begin{align*}
    \alpha_{2m-1} = \cos(y_m) + \O(m^{-1}), \hskip40pt
    \beta_{2m-1} = \sin(y_m) + \O(m^{-1}).
\end{align*}
Applying the same argument on $\psi_{2m}$ means that there exists a $\widehat{y}_m$ such that
\begin{align*}
    \alpha_{2m} = \cos(\widehat{y}_m) + \O(m^{-1}), \hskip40pt
    \beta_{2m} = \sin(\widehat{y}_m) + \O(m^{-1}).
\end{align*}
On the other hand, since $\psi_{2m-1}$ and $\psi_{2m}$ are orthogonal, we have
\begin{align*}
    0 = \langle \psi_{2m-1}, \psi_{2m} \rangle &= \alpha_{2m-1} \alpha_{2m} + \beta_{2m-1} \beta_{2m} + \O(m^{-1}) 
    = \cos(y_m - \widehat{y}_m) + \O(m^{-1}).
\end{align*}
This implies that 
\[ \widehat{y}_m = y_m + \left(k -  \frac{1}{2}\right)\pi + \O(m^{-1}) \]
which we can then apply to $\alpha_{2m}$ and $\beta_{2m}$ to yield the desired estimates, up to flipping the sign of $\psi_{2m}$.
\end{proof}

\begin{rem}
    It can be verified that the union of the two collections $\{\cos(y_m) a_m + \sin(y_m) b_m\}$ and $\{\sin(y_m) a_m - \cos(y_m) b_m\}$ 
    \begin{itemize}
        \item forms an orthonormal basis, and
        \item spans the same subspace of $L^2$ as the basis formed by the two collections $\{a_m\}$ and $\{b_m\}$.
    \end{itemize}
    In particular, we have the useful trigonometric identity
    \begin{align} \label{eq:basis_switch}
        \hspace{-5pt}\langle f, a_m \rangle a_m(x) + \langle f, b_m \rangle b_m(x) 
        &= \left( \cos(y_m) \langle f, a_m \rangle + \sin(y_m) \langle f, b_m \rangle \right) \left( \cos(y_m) a_m(x) + \sin(y_m) b_m(x) \right) \notag \\
        &\hspace{12pt} + \left( \sin(y_m) \langle f, a_m \rangle - \cos(y_m) \langle f, b_m \rangle \right) \left( \sin(y_m) a_m(x) - \cos(y_m) b_m(x) \right). 
    \end{align}
\end{rem}
    
\section{Continuous revival of the Schr\"odinger equation with piecewise $C^2$ potential}

With the asymptotics from the previous section, we are now equipped for a proof of Theorem \ref{thm:revival}. We present below a lemma that relates a solution of the periodic Schr\"odinger equation \eqref{eq:schrodinger} to the Fourier sine-cosine series of the initial datum. Theorem \ref{thm:revival} shall be a corollary to this lemma.

\begin{lem} \label{lem:revival_with_standard_basis}
    Consider a real-valued potential function $V$ that is piecewise continuously differentiable up to the second order with mean zero: $\langle V \rangle = 0$. Then a solution $u$ to the periodic Schr\"odinger equation \eqref{eq:schrodinger} has the form
    \begin{equation} \label{eq:revival_with_standard_basis}
        u(t,x) = w(t,x) + \sum_{n=1}^\infty e^{-\iu n^2t} \left[ \langle f, a_n \rangle a_n(x) + \langle f, b_n \rangle b_n(x) \right]
    \end{equation}
    where $w(t,x)$ is a continuous function for any $t > 0$.
\end{lem}

\begin{proof}
    First, consider the expansion in $L^2[0, 2\pi]$
        \[ f(x) = \sum_{n=0}^\infty c_n \psi_n(x) \roq{where } c_n = \langle f, \psi_n \rangle \]
    of the initial datum $f$ with respect to the orthonormal eigenfunction basis $\{\psi_n\}$.
    Using the eigenvalue asymptotics in \eqref{eq:eigenvalue_asymptotic_piecewise}, 
    \[ \lambda_i = m^2 + \frac{k_i}{m} \roq{for} i \in \{ 2m-1, 2m\},\]
    where $\{k_i\} \in \ell^\infty$ - that is, there exists a $K$ such that
    \[ |k_i| \leq K \]
    for all $m \geq 1$ and $i \in \{2m-1, 2m\}$. Applying this to the eigenvalue expansion of the solution to \eqref{eq:schrodinger} produces, in $L^2[0, 2\pi]$,
    \begin{align*}
        u(t,x) -  c_0 \psi_0(x)&= \sum_{n=1}^\infty c_n e^{-\iu \lambda_n t} \psi_n(x) \\
        &= \sum_{m=1}^\infty c_{2m-1} e^{-\iu (m^2 + k_{2m-1}/m) t} \psi_{2m-1}(x) + c_{2m} e^{-\iu (m^2 + k_{2m}/m) t} \psi_{2m}(x) = U_1(t,x) - U_2(t,x),
    \end{align*}
    where
    \begin{align*}
        U_1(t,x) &= \sum_{m=1}^\infty e^{-\iu m^2 t} \left( c_{2m-1} \psi_{2m-1}(x) + c_{2m} \psi_{2m}(x) \right),\\
        U_2(t,x) &= \sum_{m=1}^\infty \frac{\iu e^{-\iu m^2 t}}{m} \left( \sum_{j \in \{0, 1\}} c_{2m-j} k_{2m-j} \psi_{2m-j}(x) \int_0^t e^{-\iu k_{2m-j}s/m} ds \right).
    \end{align*}
    We shall prove that $U_2(t, \cdot) \in C([0, 2\:\pi])$ first. To do this, we first prove that the eigenfunctions are bounded. Indeed, equation \eqref{eq:eigenfunction_asymptotic_2_coeff} tells us that for all $m \geq 1$ and $i \in \{2m-1, 2m\}$, we have that
    \[ ||\psi_i||_{L^\infty[0, 2\pi]} \leq \frac{|\alpha_i| + |\beta_i|}{\sqrt{\pi}} \left[ 1 + \frac{1}{m} \int_0^{2\pi} |V(y)| \, dy \right] + ||R_i(x)||_{L^\infty[0, 2\pi]} \]
    where $||R_i||_{L^\infty} = O(m^{-2})$. Furthermore, by Lemma \ref{lem:coefficient_size}, $|\alpha_i| + |\beta_i|$ is uniformly bounded. Therefore, $||\psi_i||_{L^\infty[0, 2\pi]}$ is uniformly bounded - that is, there exists a $C_\psi$ such that
    \[ ||\psi_i||_{L^\infty[0, 2\pi]} \leq C_\psi \]
    for all $m \geq 1$ and $i \in \{2m-1, 2m\}$.
    Let
    \[ \zeta_{2m-j}(t,x) = \frac{\iu \,e^{-\iu m^2 t}}{m} c_{2m-j} k_{2m-j} \psi_{2m-j}(x) \int_0^t e^{-\iu k_{2m-j}s/m} ds. \]
    Note that for all $m \geq 1$ and $j \in \{0, 1\}$, since the sequence $\{k_{2m-j}\}$ is bounded, we have
    \begin{align*}
        |\zeta_{2m-j}(t,x)| &\leq \frac{K C_\psi  |c_{2m-j}|}{m}\, t \sup_{s \in [0, t]} \left| e^{-\iu k_{2m-j}s/m} \right| \leq \frac{C |c_{2m-j}|}{m} \,t.
    \end{align*}
    where we absorbed all constants into $C$ at the final inequality. Here, by using the Cauchy-Schwarz inequality and Parseval's identity, we get
    \begin{align*}
        \sum_{m=1}^\infty \frac{|c_{2m-1}| + |c_{2m}|}{m} &\leq \left( \sum_{m=1}^\infty \frac{1}{m^2} \right)^{1/2} \left[ \left( \sum_{m=1}^\infty |c_{2m-1}|^2 \right)^{1/2} + \left( \sum_{m=1}^\infty |c_{2m}|^2 \right)^{1/2} \right] \\
        &\leq \frac{\pi}{\sqrt{6}} \cdot 2 \left( \sum_{m=1}^\infty |c_{2m-1}|^2 + |c_{2m}|^2 \right)^{1/2} 
        \leq \frac{\pi \sqrt{6}}{3} \,||f||_{L^2[0, 2\:\pi]}.
    \end{align*}
    Therefore, for a fixed $t$, we have that
    \[ \sum_{m=1}^\infty \sum_{j \in \{0, 1\}} |\zeta_{2m-j}(t,x)| \leq C t \sum_{m=1}^\infty \sum_{j \in \{0, 1\}} \left| \frac{c_{2m-j}}{m} \right| \leq \frac{\pi \sqrt{6}}{3}\, C t ||f||_{L^2[0, 2\:\pi]} < \infty. \]
    Thus, by the Weierstrass M-test, $U_2(t, x)$ is a continuous function of $x$.
    
    We now consider $U_1(t,x)$. Using \eqref{eq:eigenfunction_asymptotic_2_coeff}, by letting
    \begin{align*}
        B_m(x) &= \sum_{j \in \{0, 1\}} c_{2m-j} \left( \alpha_{2m-j} a_m(x) + \beta_{2m-j} b_m(x) \right)
    \end{align*}
    we can write 
    \begin{align*}
        U_1(t,x) = U_3(t,x) + U_4(t,x) + U_5(t,x)
    \end{align*}
    where the component functions are
    \begin{align*}
        U_3(t,x) &= \sum_{m=1}^\infty e^{-\iu m^2 t} B_m(x), \\
        U_4(t,x) &= \sum_{m=1}^\infty \frac{e^{-\iu m^2 t}}{m} \int_0^x \sin(m(x-y)) B_m(y) V(y) \, dy, \\
        U_5(t,x) &= \sum_{m=1}^\infty e^{-\iu m^2 t} \sum_{j \in \{0, 1\}} c_{2m-j} R_{2m-j}(x).
    \end{align*}
    For $U_5(t, x)$, each individual term has the bound 
    \begin{align*}
        \left| e^{-\iu m^2 t} c_{2m-j} R_{2m-j}(x) \right| \leq ||R_{2m-j}||_{L^\infty} |\langle f, \psi_{2m-j} \rangle| \leq \frac{C ||f||_{L^2[0, 2\:\pi]}}{m^2}
    \end{align*}
    which, when applied to the Weierstrass M-test, gives us that $U_5(t,x)$ is continuous. We now turn to $U_4(t,x)$. The individual term of the series is given by
    \[ \gamma_{2m-j}(t,x) = \frac{e^{-\iu m^2 t}}{m} c_{2m-j} \int_0^x \sin(m(x-y)) [\alpha_{2m-j} a_m(y) + \beta_{2m-j} b_m(y)] V(y) \, dy. \]
    Here, with Lemma \ref{lem:coefficient_size}, for fixed $t$, there exists a constant $C$ such that
    \begin{align*}
        |\gamma_{2m-j}(t,x)| \leq \left| \frac{c_{2m-j}}{m} \right| \left(|\alpha_{2m-j}| + |\beta_{2m-j}| \right) \int_0^x |V(y)| \, dy \leq C \left| \frac{c_{2m-j}}{m} \right|
    \end{align*}
    since both coefficients $\alpha_{2m-j}$ and $\beta_{2m-j}$ are $\O(1)$ and the function $\displaystyle\int_0^x |V(y)| \, dy$ is bounded. Utilizing the same bounds obtained via the Cauchy-Schwarz inequality and Parseval's identity while evaluating $U_2(t, x)$, we get that
    \[ \sum_{m=1}^\infty \sum_{j \in \{0, 1\}} |\gamma_{2m-j}(t,x)| \leq C \sum_{m=1}^\infty \sum_{j \in \{0, 1\}} \left| \frac{c_{2m-j}}{m} \right| \leq \frac{C \pi \sqrt{6}}{3}\, ||f||_{L^2[0, 2\:\pi]} < \infty. \]
    Using the Weierstrass M-test, we deduce that $U_4(t,x)$ is a continuous function of $x$. Finally, we look at $U_3(t,x)$. From \eqref{eq:eigenfunction_asymptotic_2_coeff}, we write
    \[ c_{2m-j} = \alpha_{2m-j} \langle f, a_m \rangle + \beta_{2m-j} \langle f, b_m \rangle + \frac{1}{m} \,\langle f, \widetilde{V}_{2m-j} \rangle + \langle f, R_{2m-j} \rangle \]
    where the function $\widetilde{V}$ is given by
    \begin{align*}
        \widetilde{V}_{2m-j}(x) &= \int_0^x \sin(m(x-y)) \left( \alpha_{2m-j} a_m(y) + \beta_{2m-j} b_m(y) \right) V(y) \, dy.
    \end{align*}
    Then, we can write the function $U_3(t,x)$ as the sum
    \[ U_3(t,x) = U_6(t,x) + U_7(t,x) + U_8(t,x) \]
    where the component functions are given by
    \begin{align*}
        U_6(t,x) &= \sum_{m=1}^{\infty} e^{-\iu m^2t} \sum_{j \in \{0, 1\}} \left( \alpha_{2m-j} \langle f, a_m \rangle + \beta_{2m-j} \langle f, b_m \rangle \right) \left( \alpha_{2m-j} a_m(x) + \beta_{2m-j} b_m(x) \right),\\
        U_7(t,x) &= \sum_{m=1}^{\infty} \frac{e^{-\iu m^2t}}{m} \sum_{j \in \{0, 1\}} \langle f, \widetilde{V}_{2m-j} \rangle \left( \alpha_{2m-j} a_m(x) + \beta_{2m-j} b_m(x) \right),\\
        U_8(t,x) &= \sum_{m=1}^{\infty} e^{-\iu m^2t} \sum_{j \in \{0, 1\}} \langle f, R_{2m-j} \rangle \left( \alpha_{2m-j} a_m(x) + \beta_{2m-j} b_m(x) \right).
    \end{align*}
    For $U_8(t,x)$, thanks to the Cauchy-Schwarz inequality, we have
    \[ \langle f, R_{2m-j} \rangle = \left| \int_0^{2\:\pi} f(x) R_{2m-j}(x) \, dx \right| \leq \frac{2\:\pi C\,||f||_{L^2[0, 2\:\pi]}}{m^2} \]
    which means each individual term of $U_8(t,x)$ has the bound
    \[ \left| e^{-\iu m^2t} \langle f, R_{2m-j} \rangle \left( \alpha_{2m-j} a_m(x) + \beta_{2m-j} b_m(x) \right) \right| \leq \frac{C ||f||_{L^2[0, 2\:\pi]}}{m^2} \]
    for some constant $C$, as the coefficients $\alpha_{2m-j}$ and $\beta_{2m-j}$ are $\O(1)$. Applying the Weierstrass M-test here shows that $U_8(t,x)$ is a continuous function. For $U_7(t,x)$, note that
    \begin{align*}
        \widetilde{V}_{2m-j}(x) &= \sin(mx) \int_0^x \cos(my) \left( \alpha_{2m-j} a_m(y) + \beta_{2m-j} b_m(y) \right) V(y) \, dy \\
        &\hspace{55pt} - \cos(mx) \int_0^x \sin(my) \left( \alpha_{2m-j} a_m(y) + \beta_{2m-j} b_m(y) \right) V(y) \, dy
    \end{align*}
    where the first term can be written as
    \begin{align*}
         \sin(mx) \int_0^x \cos(&my) \left( \alpha_{2m-j} a_m(y) + \beta_{2m-j} b_m(y) \right) V(y) \, dy \\
        &= \frac{\sin(mx)}{2\sqrt{\pi}} \int_0^x \left( \alpha_{2m-j} (1 + \cos(2my)) + \beta_{2m-j} \sin(2my)\right) V(y) \, dy \\
        &= \frac{\alpha_{2m-j}}{2} b_m(x) \int_0^x V(y) \, dy + \frac{b_m(x)}{2} \int_0^x \left( \alpha_{2m-j}\cos(2my) + \beta_{2m-j} \sin(2my)\right) V(y) \, dy
    \end{align*}
    and the second term can be written as
    \begin{align*}
        \cos(mx) \int_0^x \sin(&my) \left( \alpha_{2m-j} a_m(y) + \beta_{2m-j} b_m(y) \right) V(y) \, dy \\
        &= \frac{\cos(mx)}{2\sqrt{\pi}} \int_0^x \left( \alpha_{2m-j} \sin(2my) + \beta_{2m-j} (1-\cos(2my)) \right) V(y) \, dy \\
        &= \frac{\beta_{2m-j}}{2} a_m(x) \int_0^x V(y) \, dy + \frac{a_m(x)}{2} \int_0^x \left( \alpha_{2m-j} \sin(2my) - \beta_{2m-j}\cos(2my) \right) V(y) \, dy.
    \end{align*}
    Now, using integration by parts, with $N_0$ such that $x_{N_0} < x < x_{N_0 + 1}$ --- that is, the $N_0$\textsuperscript{th} discontinuity point is the greatest discontinuity point of $V(x)$ smaller than $x$ --- we have
    \begin{align*}
        \int_0^x \sin(2my) V(y) \, dy &= - \frac{1}{2m} \sum_{n=0}^{N_0 - 1} \left( \cos(2mx_{n+1}) \lim_{x \to x_{n+1}^-} V(x) - \cos(2mx_n) \lim_{x \to x_n^+} V(x) \right) \\
        &\hspace{11pt} - \frac{1}{2m} \left( \cos(2mx) V(x) - \cos(2mx_{N_0}) \lim_{x \to x_{N_0}^+} V(x) \right) + \frac{1}{2m} \int_0^{x} \cos(2my) V'(y) \, dy 
    \end{align*}
    which grows like $\O(m^{-1})$. This also applies to
    \begin{align*}
        \int_0^x \cos(2my) V(y) \, dy &= \frac{1}{2m} \sum_{n=0}^{N_0 - 1} \left( \sin(2mx_{n+1}) \lim_{x \to x_{n+1}^-} V(x) - \sin(2mx_n) \lim_{x \to x_n^+} V(x) \right) \\
        &\hspace{11pt} + \frac{1}{2m} \left( \sin(2mx) V(x) - \sin(2mx_{N_0}) \lim_{x \to x_{N_0}^+} V(x) \right) - \frac{1}{2m} \int_0^{x} \sin(2my) V'(y) \, dy.
    \end{align*}
    Then, each individual term of $U_7(t,x)$ can be split into three terms
    \[ \frac{e^{-\iu m^2t}}{m} \langle f, \widetilde{V}_{2m-j} \rangle \left( \alpha_{2m-j} a_m(x) + \beta_{2m-j} b_m(x) \right) = \rho_{1, 2m-j}(t,x) + \rho_{2, 2m-j}(t,x) + \rho_{3, 2m-j}(t,x) \]
    where
    \begin{align*}
        \rho_{1, 2m-j}(t,x) &= \frac{\alpha_{2m-j} e^{-\iu m^2t}}{2m} \left \langle f V_0, b_m \right \rangle \left( \alpha_{2m-j} a_m(x) + \beta_{2m-j} b_m(x) \right) \\
        \rho_{2, 2m-j}(t,x) &= -\frac{\beta_{2m-j} e^{-\iu m^2t}}{2m} \left \langle f V_0, a_m \right \rangle \left( \alpha_{2m-j} a_m(x) + \beta_{2m-j} b_m(x) \right) \\
        \rho_{3, 2m-j}(t,x) &= e^{-\iu m^2t} \left \langle f, R^*_{2m-j} \right \rangle \left( \alpha_{2m-j} a_m(x) + \beta_{2m-j} b_m(x) \right)
    \end{align*}
    with
    \[ V_0(x) = \int_0^{x} V(y) \, dy \roq{and} ||R^*_{2m-j}||_{L^\infty} = \O(m^{-2}).\]
    As for the third sub-term, employing the same argument as with $U_8(t,x)$, we have the bound
    \[ |\rho_{3, 2m-j}(t,x)| \leq \frac{C\, ||f||_{L^2[0, 2\:\pi]}}{m^2} \]
    which means that the Weierstrass M-test gives us that the sum across all of these sub-terms is a continuous function. For the first sub-term, with $\alpha_{2m-j}$ and $\beta_{2m-j}$ both being $\O(1)$, plus that $\{b_m\}$ is the orthonormal sine basis and $fV_0 \in L^2[0, 2\:\pi]$, we have
    \[ \sum_{m=1}^\infty \sum_{j \in \{0, 1\}} \left| \rho_{1, 2m-j}(t,x) \right| \leq \sum_{m=1}^\infty C \left| \frac{\left \langle f V_0, b_m \right \rangle}{m} \right| \leq \frac{C\: \pi}{\sqrt{6}}\, ||fV_0||_{L^2[0, 2\:\pi]}. \]
    This ensures that the sum across all of these sub-terms is continuous due to the Weierstrass M-test. With the same method applied to the second sub-term, we conclude that
    \[ U_7(t,x) = \sum_{m=1}^\infty \sum_{j \in \{0, 1\}} \bigl[\,\rho_{1, 2m-j}(t,x) + \rho_{2, 2m-j}(t,x) + \rho_{3, 2m-j}(t,x)\,\bigr]\]
    is a continuous function. This leaves us with $U_6(t,x)$. Using the asymptotics from Lemma \ref{lem:coeff_asymptotic}, we can split it into two parts
    \[ U_6(t,x) = U_9(t,x) + U_{10}(t,x) \]
    where the component functions are 
    \begin{align*}
        U_9(t,x) &= \sum_{m=1}^\infty e^{-\iu m^2t} (\langle f, \cos(y_m) a_m + \sin(y_m) b_m \rangle \left( \cos(y_m) a_m(x) + \sin(y_m) b_m(x) \right) \\
        &\hspace{100pt} + \langle f, \sin(y_m) a_m - \cos(y_m) b_m \rangle \left( \sin(y_m) a_m(x) - \cos(y_m) b_m(x) \right)), \\
        U_{10}(t,x) &= \sum_{m=1}^\infty e^{-\iu m^2t} \sum_{j \in \{0, 1\}} S_{1, 2m-j} \langle f, a_m \rangle a_m(x) + S_{2, 2m-j} \langle f, a_m \rangle b_m(x) \\
        &\hspace{100pt} + S_{3, 2m-j} \langle f, b_m \rangle a_m(x) + S_{4, 2m-j} \langle f, b_m \rangle b_m(x)
    \end{align*}
    with $|S_{i, 2m-j}|$ being $\O(m^{-1})$ for $1 \leq i \leq 4$. The identity \eqref{eq:basis_switch} immediately identifies $U_9(t,x)$ as the revival component of the solution, while for $U_{10}(t,x)$, we can apply the same technique on the sum across each of the four sub-terms as we did for $U_7(t,x)$, which will lead to convergence to a continuous function. For example, the sum across the first sub-term has the bound
    \[ \sum_{m=1}^\infty \sum_{j \in \{0, 1\}} |e^{-\iu m^2t} S_{1, 2m-j} \langle f, a_m \rangle a_m| \leq \sum_{m=1}^\infty C \left| \frac{\left \langle f, a_m \right \rangle}{m} \right| \leq \frac{C \pi}{\sqrt{6}} ||f||_{L^2[0, 2\:\pi]}, \]
    which proves that $U_{10}(t,x)$ is continuous. Assembling all the components gives us that
    \[ u(t,x) = c_0 \psi_0(x) + U_1(t,x) - U_9(t,x) - U_2(t,x) + \sum_{n=1}^\infty e^{-\iu n^2t} \left[ \langle f, a_n \rangle a_n(x) + \langle f, b_n \rangle b_n(x) \right]. \]
    Setting $w(t,x) = c_0 \psi_0(x) + U_1(t,x) - U_9(t,x) - U_2(t,x)$ then gives us the desired identity.
\end{proof}

The proof of Theorem \ref{thm:revival} now follows from this lemma and a rearrangement with regard to the norm of the potential function.

\begin{proof}[Proof of Theorem \ref{thm:revival}]
    First, for $V = 0$, the solution to \eqref{eq:schrodinger} is given by
    \[ u_0(t,x) = \frac{\langle f^*, 1 \rangle}{2\pi} + \sum_{n=1}^\infty e^{-\iu n^2t} \left[\, \langle f^*, a_n \rangle a_n(x) + \langle f^*, b_n \rangle b_n(x) \,\right] = \frac{1}{2\:\pi} \sum_{n=-\infty}^\infty e^{-\iu n^2t} \langle f^*, e^{\iu n(\cdot)} \rangle e^{\iu nx} . \]
    in $L^2[0, 2\pi]$. Then, at a rational time $t_0 = 2\:\pi q/r$, observe that
    \[ e^{\iu n^2 t_0} = e^{\iu m^2 t_0} \]
    for any $n \equiv m \bmod r$. Thus, the solution has the form
    \[ u_0(t_0,x) = \frac{1}{2\:\pi} \sum_{m=0}^{r-1} e^{-2\:\pi \iu m^2q/r} \sum_{n \equiv m \bmod r} \langle f^*, e^{\iu n(\cdot)} \rangle e^{\iu nx}. \]
    We shall focus on the inner summation. Since
    \begin{align*}
        \sum_{k=0}^{r-1} e^{2\:\pi \iu (m-n)k/r} = \begin{cases}
            r \qquad &\roq{if} n \equiv m \bmod r, \\
            0 &\roq{otherwise,}
        \end{cases}
    \end{align*}
    we have, in $L^2[0, 2\pi]$, 
    \begin{align*}
        \sum_{n \equiv m \bmod r} \langle f^*, e^{\iu n(\cdot)} \rangle e^{\iu nx} &= \frac{1}{r} \sum_{k=0}^{r-1} e^{2\:\pi \iu mk/r} \sum_{n \in \ZZ} e^{-2\:\pi \iu nk/r} \langle f^*, e^{\iu n(\cdot)} \rangle e^{\iu nx} \\
        &= \frac{1}{r} \sum_{k=0}^{r-1} e^{2\:\pi \iu mk/r} \sum_{n \in \ZZ} \left \langle f^*\left(\cdot - \frac{2\:\pi k}{r} \right), e^{\iu n(\cdot)} \right \rangle e^{\iu nx} 
        = \frac{2\:\pi}{r} \sum_{k=0}^{r-1} e^{2\:\pi \iu mk/r} f^*\left(x - \frac{2\:\pi k}{r} \right).
    \end{align*}
    Substituting this in the full expression above for the solution gives
    \[ u_0(t_0,x) = \frac{1}{r} \sum^{r-1}_{k, m=0} e^{2\:\pi \iu \left(mk/r- m^2q/r \right)} f^*\left(x - 2\:\pi \frac{k}{r} \right). \]
    By \eqref{eq:revival_with_standard_basis}, for any potential $V$ satisfying $\langle V \rangle = 0$,
    \begin{align*}  
        u(t_0,x) &= w(t_0,x) + \left[ u_0(t_0,x) - \frac{\langle f^*, 1 \rangle}{2\pi} \right]
        = \left[w(t_0,x) - \frac{\langle f^*, 1 \rangle}{2\pi} \right] + \frac{1}{r} \sum^{r-1}_{k, m=0} e^{2\:\pi \iu \left(mk/r- m^2q/r \right)} f^*\left(x - 2\:\pi \frac{k}{r} \right).
    \end{align*}
    Finally, notice that for a general potential $V$, the solution to \eqref{eq:schrodinger} has the form
    \[ u(t,x) = e^{-\iu \langle V \rangle t} v(t,x) \]
    where $v$ is the solution to
    \begin{align*}
        \begin{cases}
            \iu \, v_t = - v_{xx} + \left[ V - \langle V \rangle \right] v, \\[1pt]
            v(0, x) = f(x),
        \end{cases}
    \end{align*}
    which has a potential that satisfies the mean-zero condition. Applying the derived result above, we have that
    \[ u(t_0,x) = e^{-2\:\pi \iu \langle V \rangle q/r} \left( w(t_0,x) - \frac{\langle f^*, 1 \rangle}{2\pi} \right) + \frac{1}{r} e^{-2\:\pi \iu \langle V \rangle q/r} \sum^{r-1}_{k, m=0} e^{2\:\pi \iu \left(mk/r- m^2q/r \right)} f^*\left(x - 2\:\pi \frac{k}{r} \right). \]
    Since the first term is a continuous function, the proof is complete.
\end{proof}

\section{Numerical simulations}

In this section, we wish to show various simulations of the solution to the Schr\"odinger equation. Consider the case where the potential  is the following $2\:\pi$-periodic piecewise constant function:
\begin{align*}
    V(x) = \begin{cases}
        0 \qquad &\roq{for} x \in [0, \pi / 2), \\
        9 &\roq{for} x \in [\pi / 2, 2\:\pi),
    \end{cases}
\end{align*}
and the initial datum is the following $2\:\pi$-periodic function:
\begin{align*}
    f(x) = \begin{cases}
        - x/2\:\pi \qquad &\roq{for} x \in [0, \pi), \\
        1 - x/2\:\pi &\roq{for} x \in [\pi, 2\:\pi).
    \end{cases}
\end{align*}

Here, our plots were produced using Julia employing these assumptions:
\begin{itemize}
    \setlength{\itemsep}{3pt} 
    \item The number of eigenvalues used is 1000;
    \item The resolution of the plot is $\Delta x = 0.0005\pi$.
\end{itemize}

As a dispersive relation, the solution of the Schr\"odinger equation is piecewise smooth for times rational relative to the length of the interval. This can be seen in the plots of the real part of the solution below.

\begin{figure}[H]
    \begin{tabularx}{\textwidth}{>{\centering\arraybackslash}X >{\centering\arraybackslash}X >{\centering\arraybackslash}X >{\centering\arraybackslash}X} 
        \includegraphics[width = 4cm]{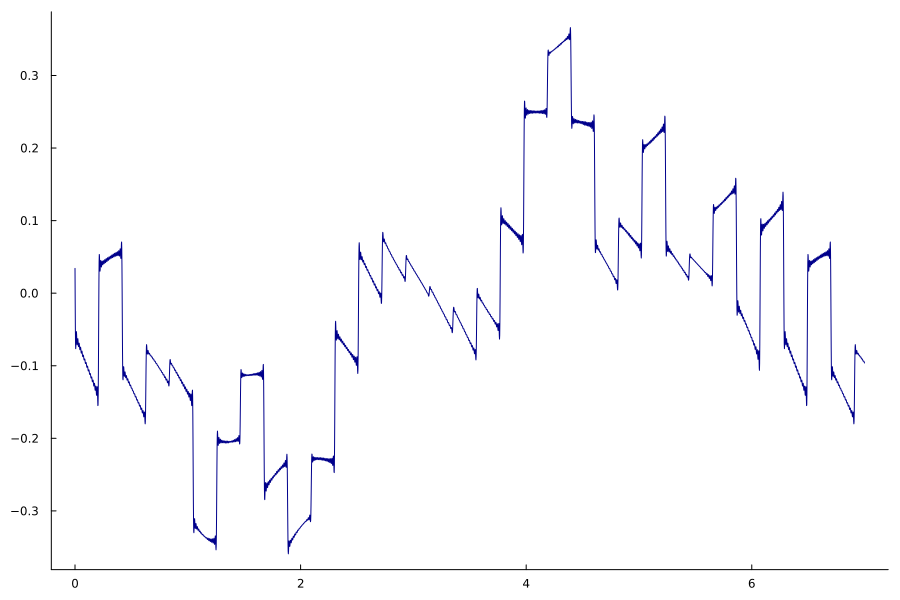}
        $t = \frac{1}{30} \pi$
        & \includegraphics[width = 4cm]{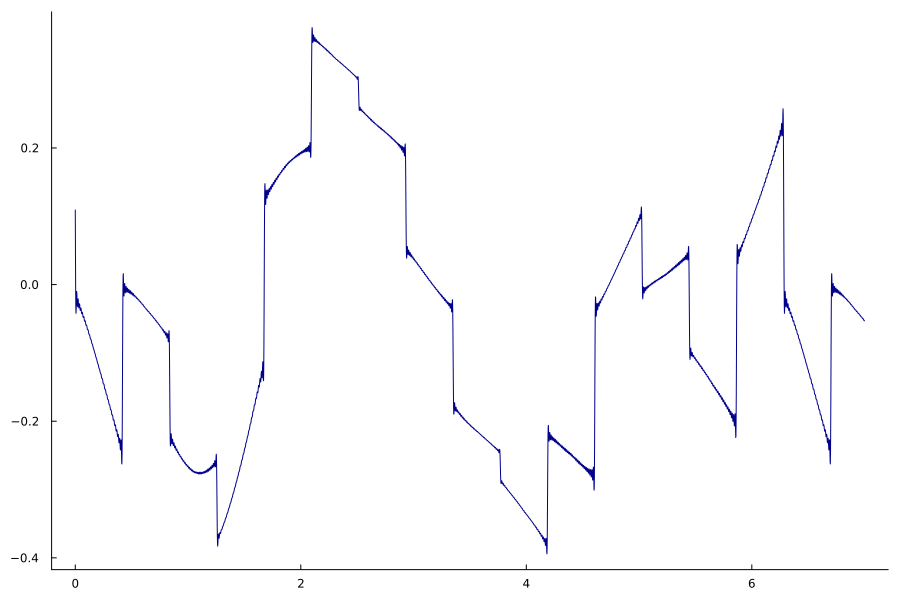}
        $t = \frac{1}{15} \pi$
        & \includegraphics[width = 4cm]{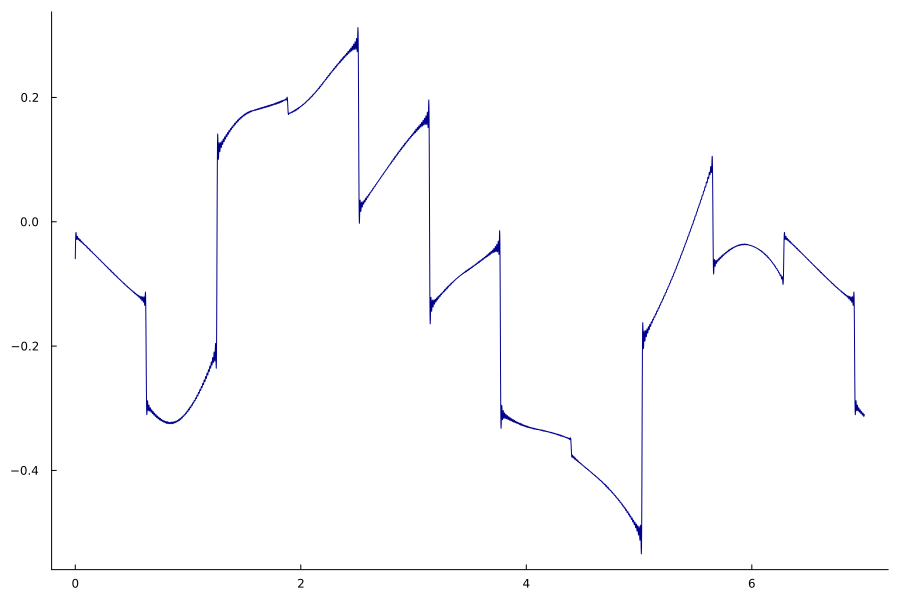}
        $t = \frac{1}{10} \pi$
        & \includegraphics[width = 4cm]{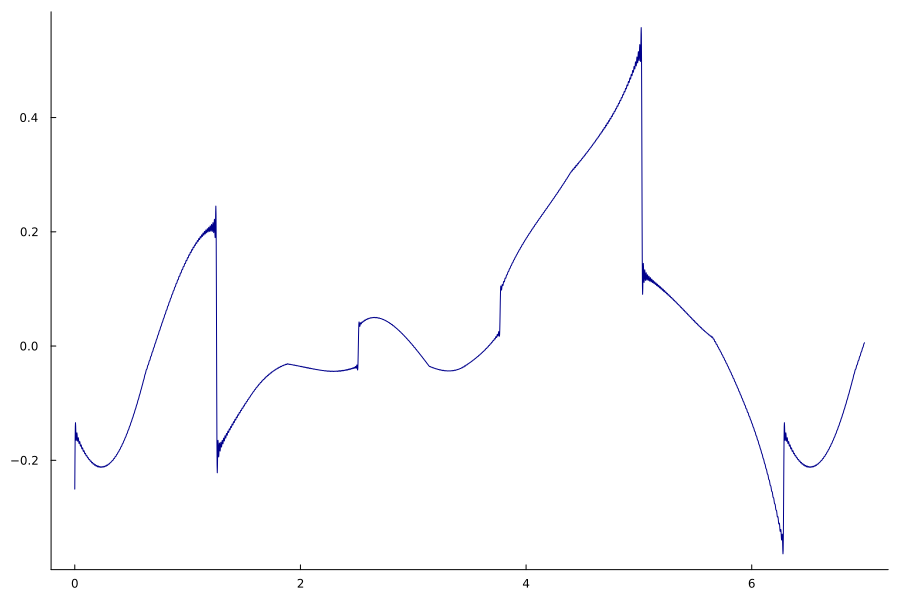}
        $t = \frac{1}{5} \pi$
    \end{tabularx}
    \caption{The real part of the solution.}
\end{figure}

We then consider the plots below of the same rational times, where we have the real part of the revival component and the remainder $w$ on the top and bottom row, respectively. 

\begin{figure}[H]
    \begin{tabularx}{\textwidth}{>{\centering\arraybackslash}X >{\centering\arraybackslash}X >{\centering\arraybackslash}X >{\centering\arraybackslash}X} 
        \includegraphics[width = 4cm]{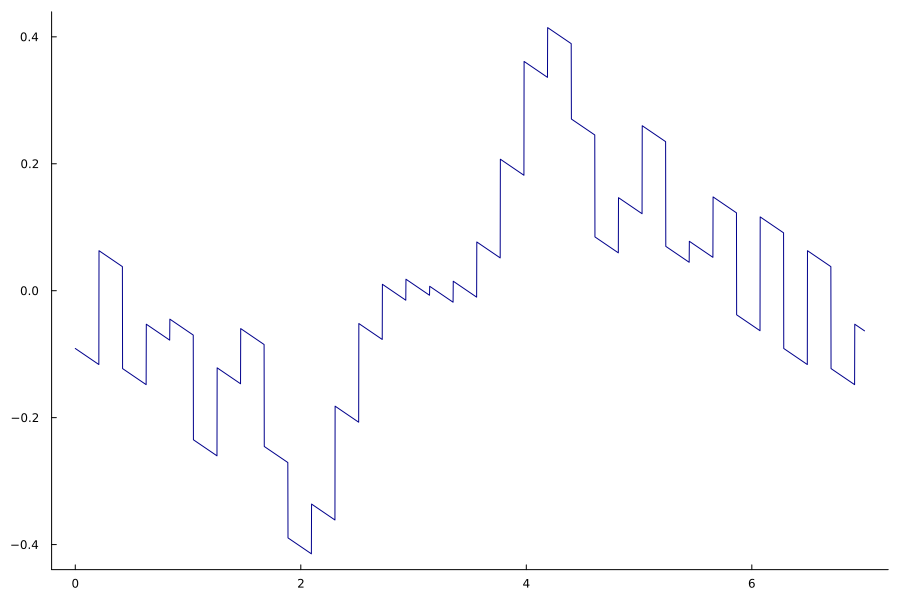}
        $t = \frac{1}{30} \pi$
        & \includegraphics[width = 4cm]{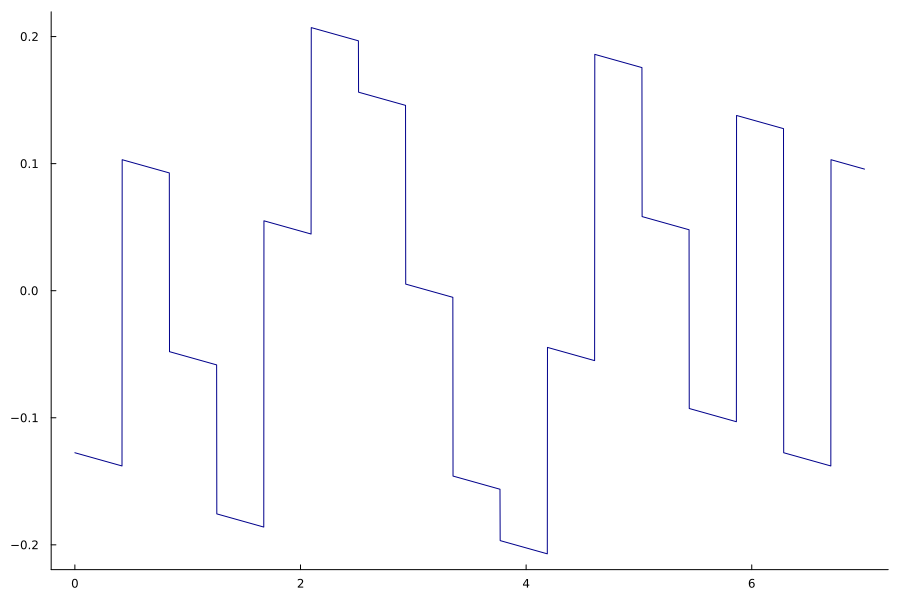}
        $t = \frac{1}{15} \pi$
        & \includegraphics[width = 4cm]{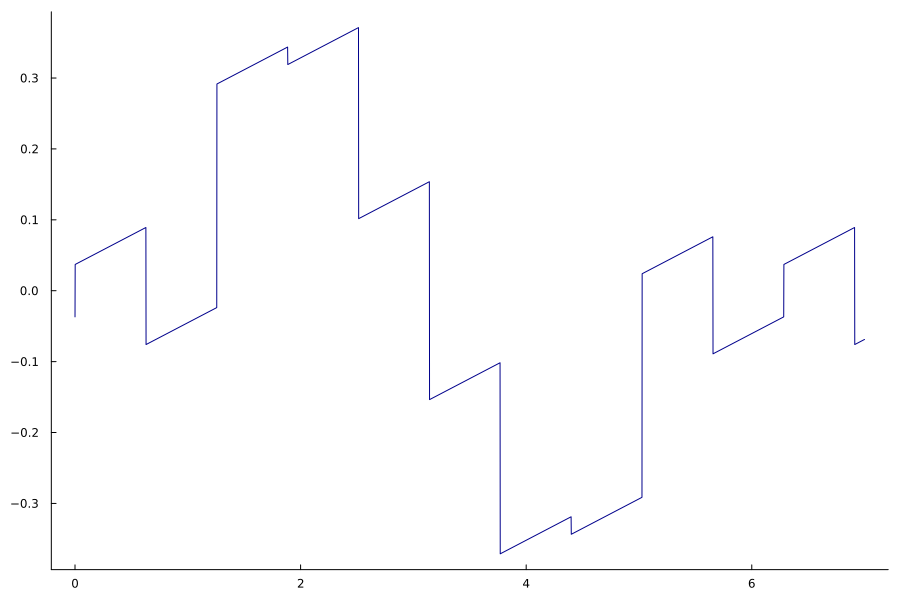}
        $t = \frac{1}{10} \pi$
        & \includegraphics[width = 4cm]{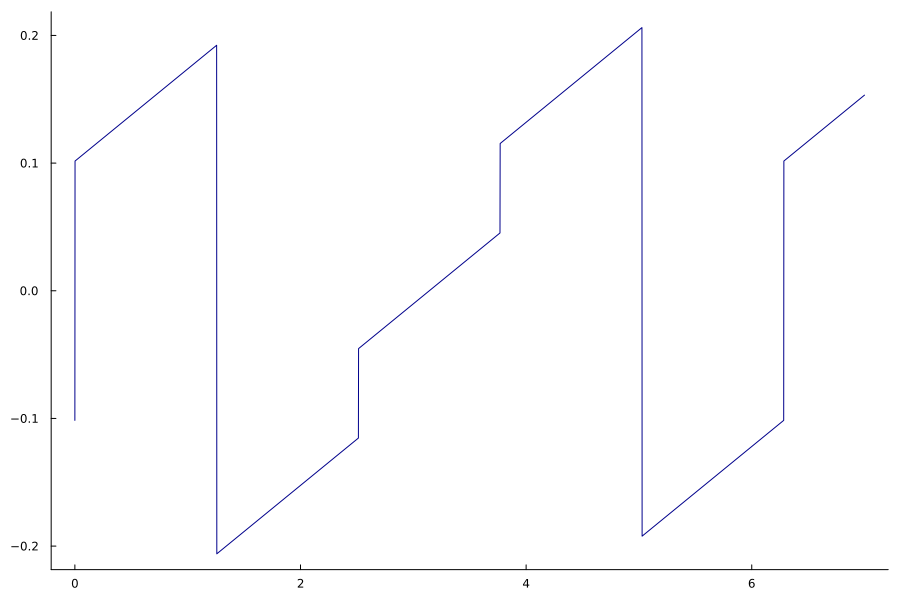}
        $t = \frac{1}{5} \pi$
    \end{tabularx}
\end{figure}

\begin{figure}[H]
    \begin{tabularx}{\textwidth}{>{\centering\arraybackslash}X >{\centering\arraybackslash}X >{\centering\arraybackslash}X >{\centering\arraybackslash}X} 
        \includegraphics[width = 4cm]{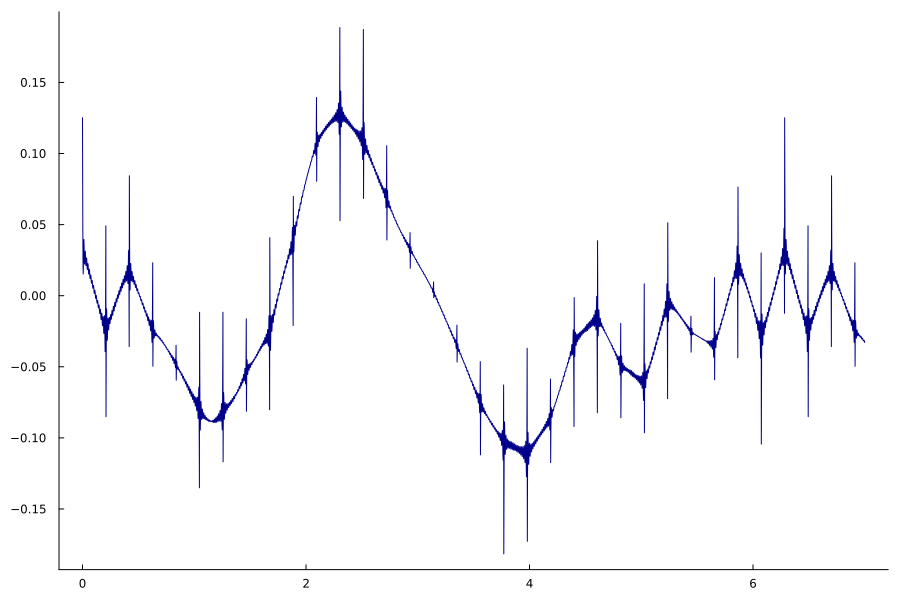}
        $t = \frac{1}{30} \pi$
        & \includegraphics[width = 4cm]{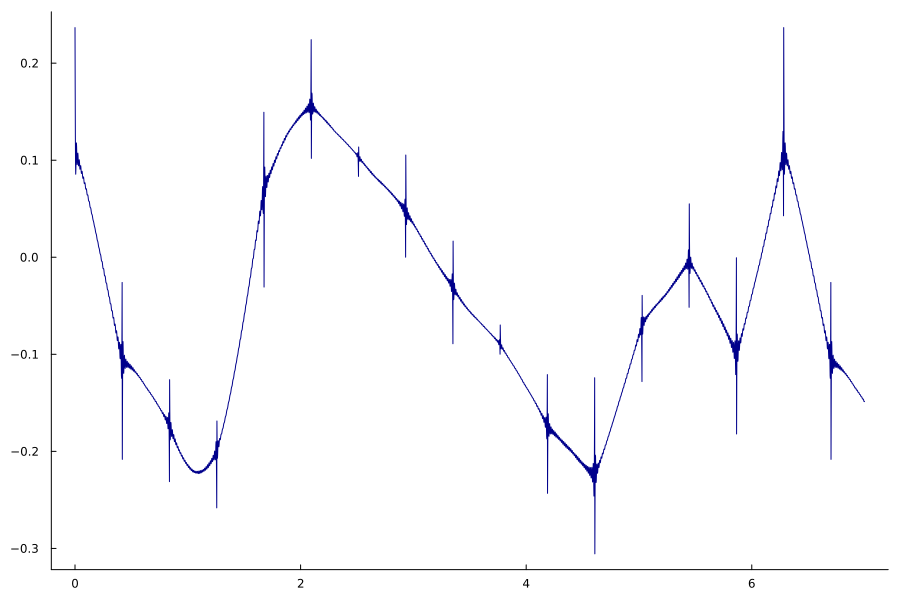}
        $t = \frac{1}{15} \pi$
        & \includegraphics[width = 4cm]{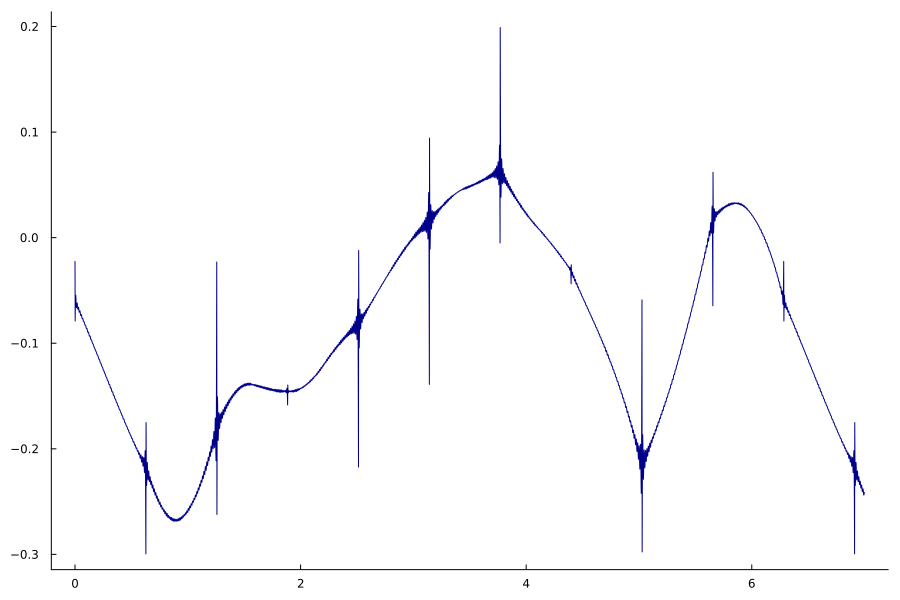}
        $t = \frac{1}{10} \pi$ 
        & \includegraphics[width = 4cm]{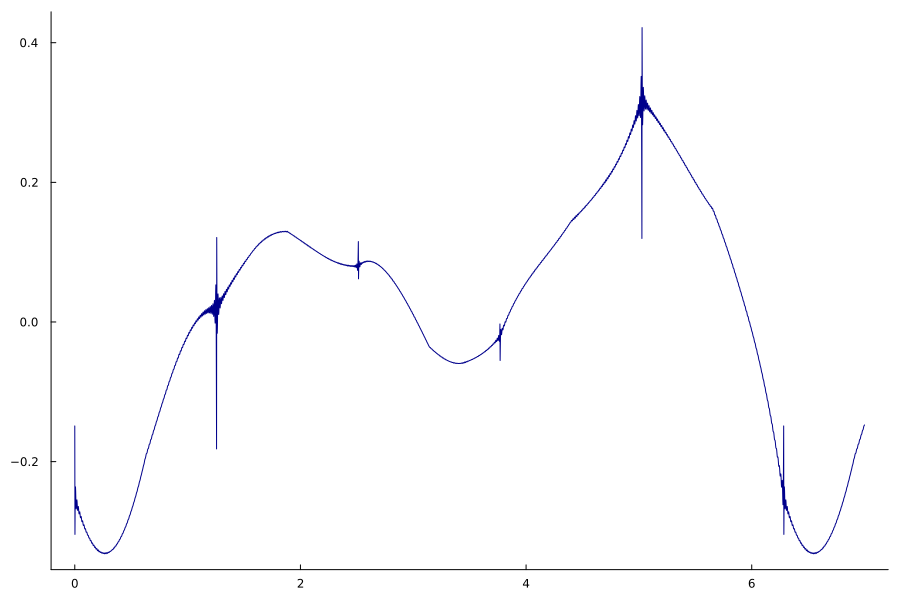}
        $t = \frac{1}{5} \pi$
    \end{tabularx}

    \caption{The real part of the revival component and the remainder.}
\end{figure}

Apart from the bars at the jump discontinuities caused by the Gibbs phenomenon, it can be observed that the real part of $w$ is continuous across all provided rational times. Similarly, we demonstrate below the continuity for the imaginary part of $w$, which is plotted in the bottom row and is the difference between the imaginary part of the solution (top row) and the imaginary part of the revival component (middle row).

\begin{figure}[H]
    \begin{tabularx}{\textwidth}{>{\centering\arraybackslash}X >{\centering\arraybackslash}X >{\centering\arraybackslash}X >{\centering\arraybackslash}X} 
        \includegraphics[width = 4cm]{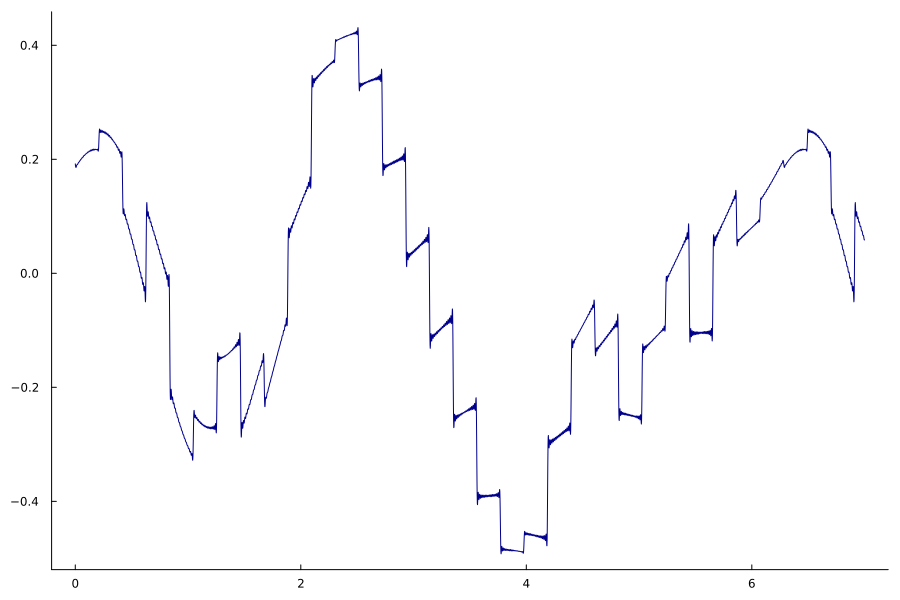}
        $t = \frac{1}{30} \pi$
        & \includegraphics[width = 4cm]{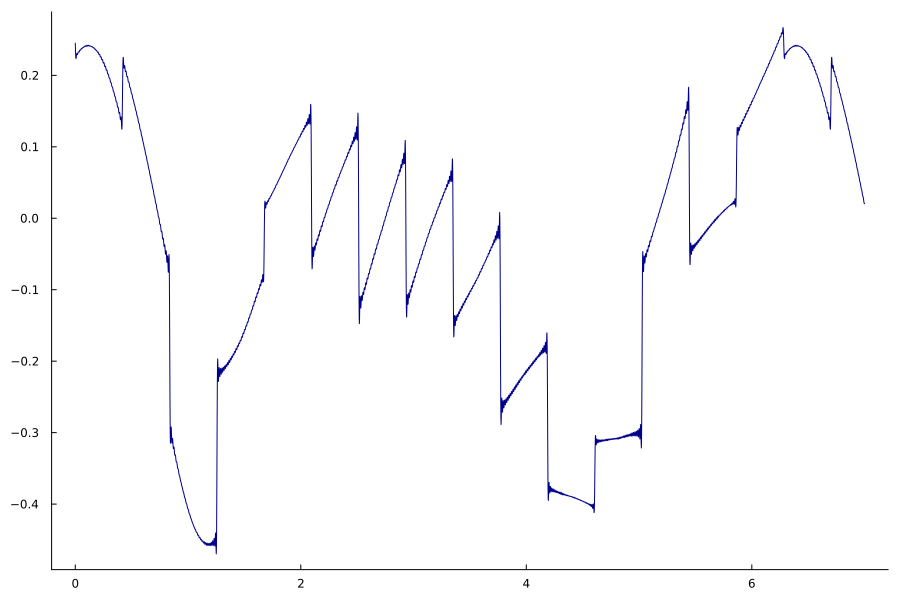}
        $t = \frac{1}{15} \pi$
        & \includegraphics[width = 4cm]{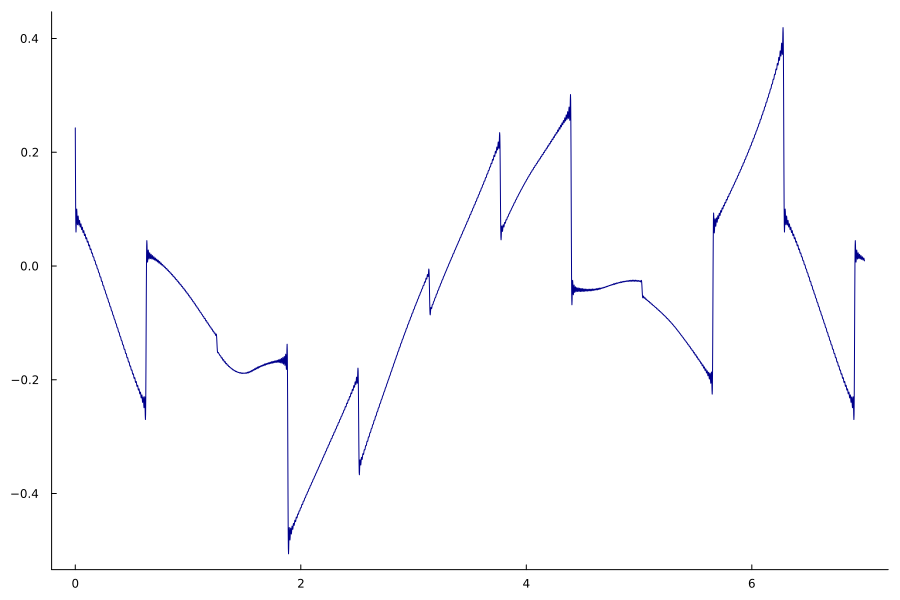}
        $t = \frac{1}{10} \pi$
        & \includegraphics[width = 4cm]{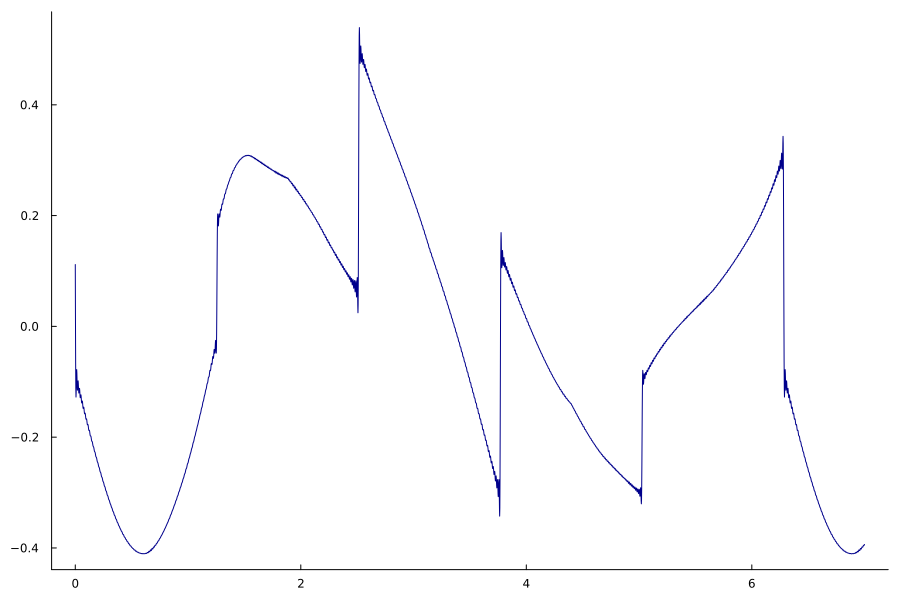}
        $t = \frac{1}{5} \pi$
    \end{tabularx}
\end{figure}

\begin{figure}[H]
    \begin{tabularx}{\textwidth}{>{\centering\arraybackslash}X >{\centering\arraybackslash}X >{\centering\arraybackslash}X >{\centering\arraybackslash}X} 
        \includegraphics[width = 4cm]{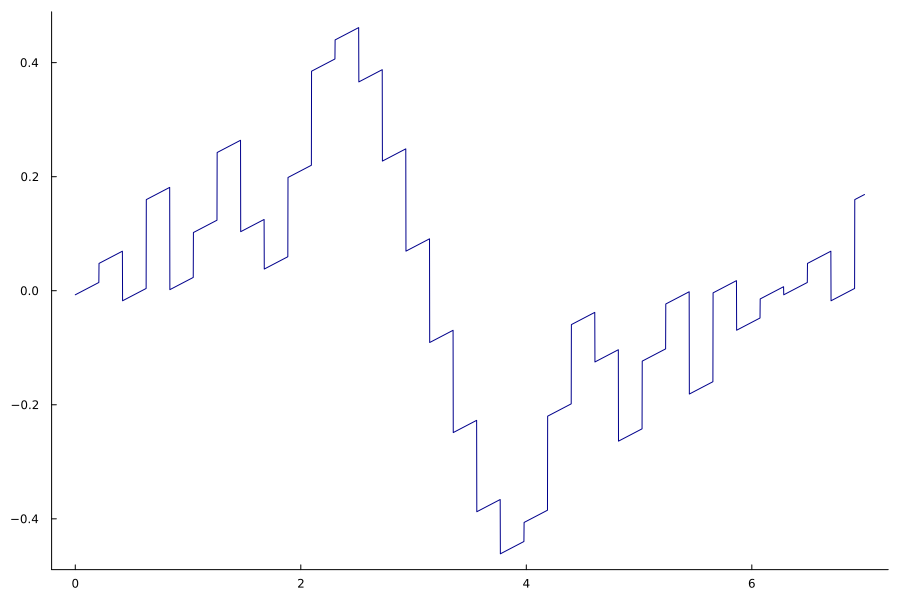}
        $t = \frac{1}{30} \pi$
        & \includegraphics[width = 4cm]{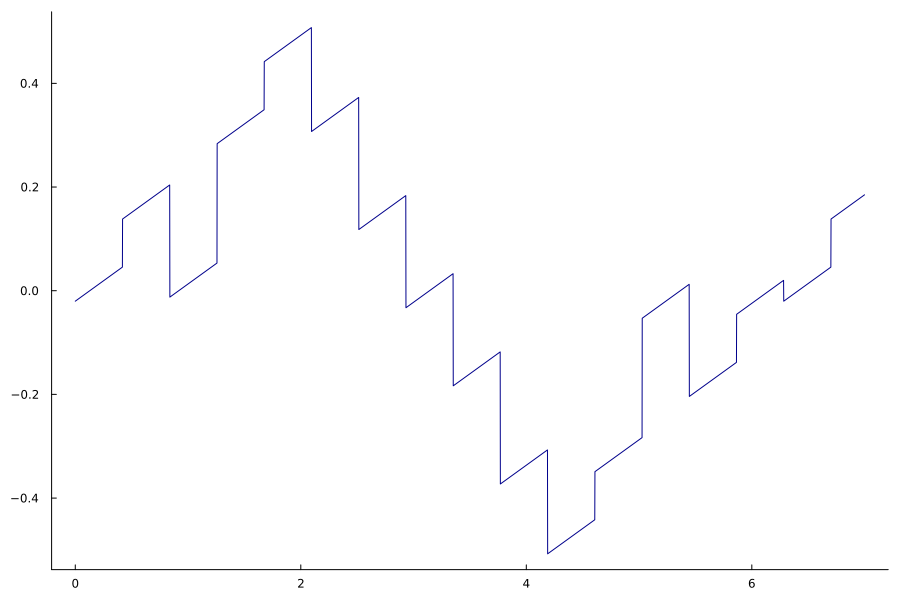}
        $t = \frac{1}{15} \pi$
        & \includegraphics[width = 4cm]{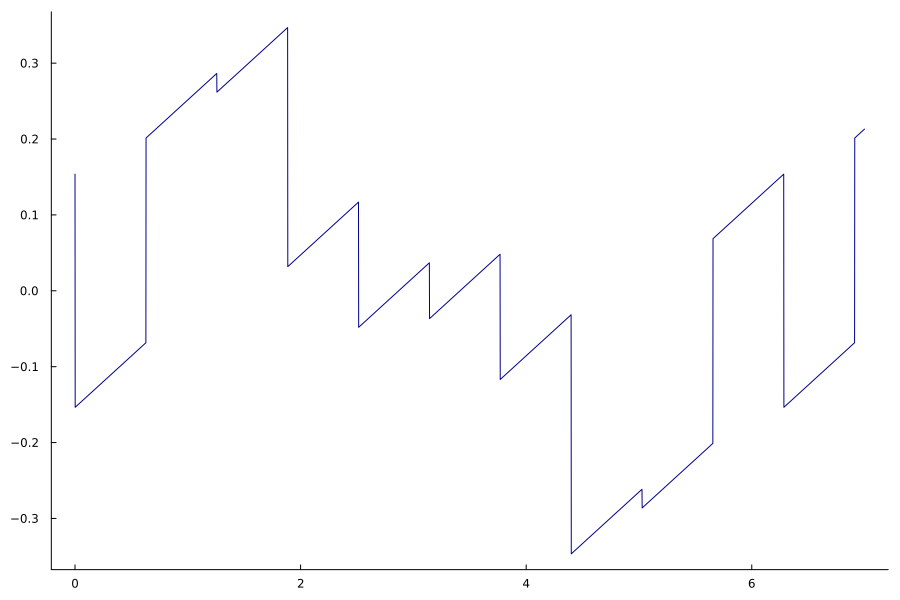}
        $t = \frac{1}{10} \pi$
        & \includegraphics[width = 4cm]{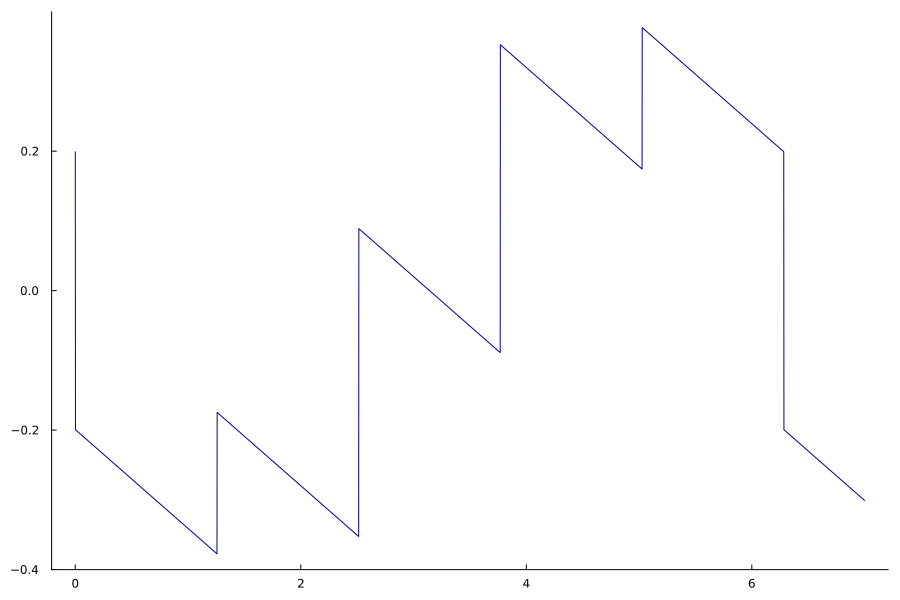}
        $t = \frac{1}{5} \pi$
    \end{tabularx}
\end{figure}

\begin{figure}[H]
    \begin{tabularx}{\textwidth}{>{\centering\arraybackslash}X >{\centering\arraybackslash}X >{\centering\arraybackslash}X >{\centering\arraybackslash}X} 
        \includegraphics[width = 4cm]{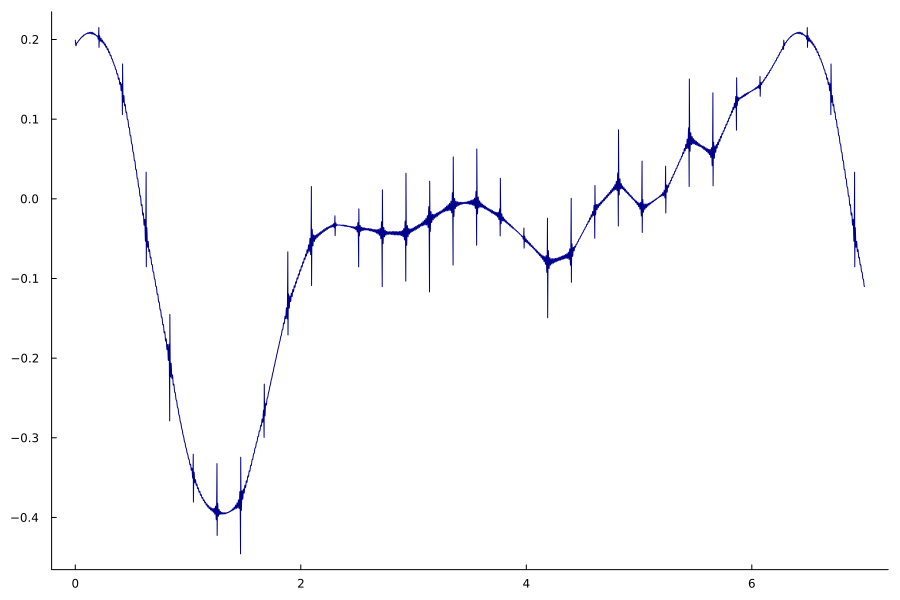}
        $t = \frac{1}{30} \pi$
        & \includegraphics[width = 4cm]{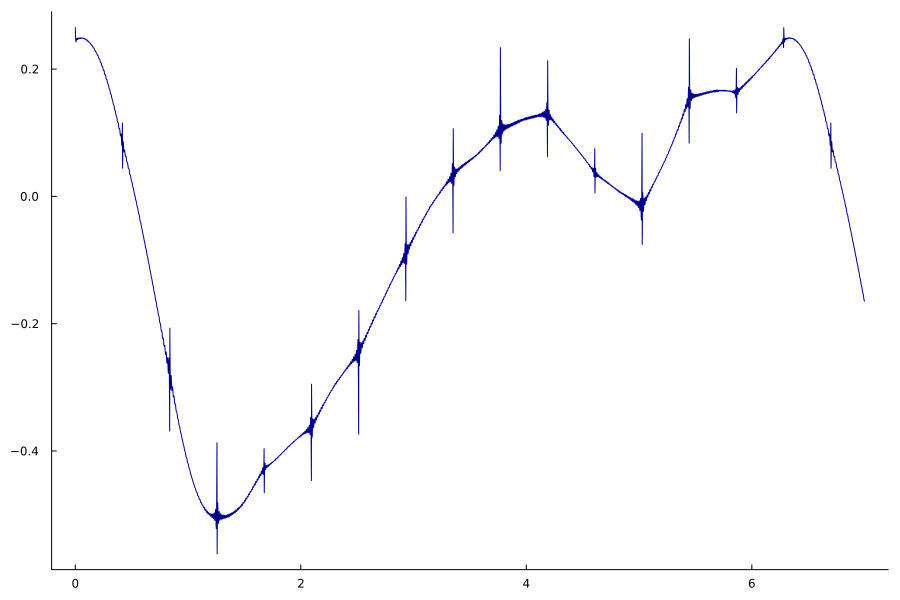}
        $t = \frac{1}{15} \pi$
        & \includegraphics[width = 4cm]{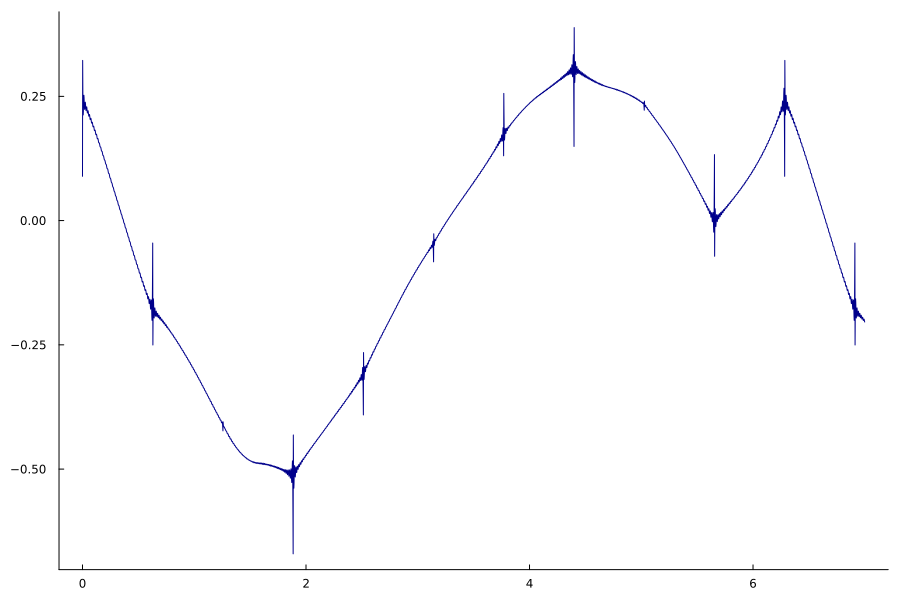}
        $t = \frac{1}{10} \pi$ 
        & \includegraphics[width = 4cm]{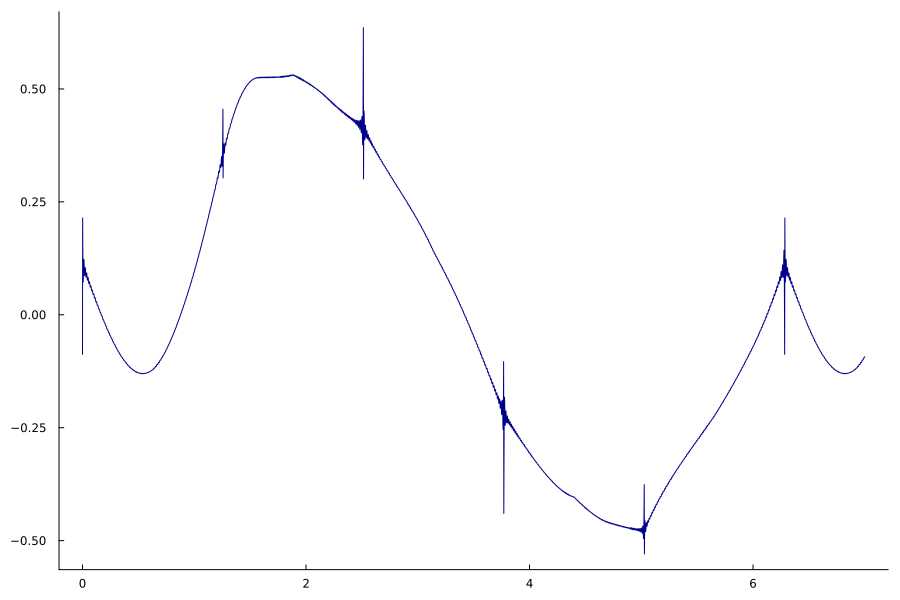}
        $t = \frac{1}{5} \pi$
    \end{tabularx}
    \caption{The imaginary part of the solution, the revival component, and the remainder.}
\end{figure}

This illustrates the consistency of Theorem \ref{thm:revival} in this particular example.

\section*{Acknowledgements}

The authors thank Jing Kang, Changzheng Qu, and Zihan Yin for their input on these results. We also thank the referee for carefully reading the paper, constructive comments, and for noticing several typos and pointing out a mathematical error in the original version.

\newpage

\bibliographystyle{plain}
\bibliography{references}

@article{chen_olver_2014,
    author = {Chen, Gong and Olver, Peter J.},
    title = {Numerical {Simulation} of {Nonlinear} {Dispersive} {Quantization}},
    journal = {Discrete and Continuous Dynamical Systems},
    year = {2014},
    volume = {34},
    number = {3},
    pages = {991--1008},
}

@article{olver_2010,
    author = {Olver, Peter J.},
    title = {Dispersive {Quantization}},
    journal = {The American Mathematical Monthly},
    year = {2010},
    volume = {117},
    number = {7},
    pages = {599--610},
}

@article{boulton_farmakis_pelloni_2024,
    author={Boulton, Lyonell and Farmakis, George and Pelloni, Beatrice},
    title = {The {Phenomenon} of {Revivals} on {Complex} {Potential} {Schr\"odinger's} {Equation}},
    journal = {Zeitschrift f{\"u}r Analysis und ihre Anwendungen},
    year = {2024},
    volume = {43},
    number = {3-4},
    pages = {401--416},
}

@book{eastham_1973,
    author={Eastham, M. S. P.},
    title={{The Spectral Theory of Periodic Differential Equations}},
    series={Texts in Mathematics},
    year={1973},
    publisher={Scottish Academic Press, London},
}

@incollection{simon_2005,
    author = {Simon, Barry},
    editor = {Amrein, Werner O. and Hinz, Andreas M. and Pearson, David P.},
    title = {Sturm {Oscillation} and {Comparison} {Theorems}},
    booktitle = {Sturm--{Liouville} {Theory}: {Past} and {Present}},
    year = {2005}, 
    publisher = {Birkh{\"a}user Basel}, 
    pages = {29--43},
}

@book{eastham_ode_1970,
    author = {Eastham, M. S. P.},
    title = {Theory of {Ordinary} {Differential} {Equations}},
	series = {New University Mathematics Series},
	year = {1970},
    publisher = {Van Nostrand Reinhold},
}

@article{prufer_1926,
    author = {Pr{\"u}fer, H.},
    title = {{Neue Herleitung der Sturm-Liouvilleschen Reihenentwicklung stetiger Funktionen}},
    journal = {Mathematische Annalen},
    year = {1926},
    volume = {95},
    pages = {499--518},
}

@article{berry_marzoli_schleich_2001,
    author = {Berry, Michael and Marzoli, Irene and Schleich, Wolfgang},
    title = {{Quantum Carpets, Carpets of Light}},
    journal = {Physics World},
    year = {2001},
    volume = {14},
    number = {6},
    pages = {39--44},
}

@article{berry_1996,
    author = {Berry, M. V.},
    title = {{Quantum Fractals in Boxes}},
    journal = {Journal of Physics A: Mathematical and General},
    year = {1996},
    volume = {29},
    number = {20},
    pages = {6617--6629},
}

@article{cho_2024,
    author = {Gunwoo Cho and Jimyeong Kim and Ihyeok Seo and Seongyeon Kim and Yehyun Kwon},
    title = {{Quantum Revivals and Fractality for the Schr\"{o}dinger Equation}},
    journal = {Reports on Mathematical Physics},
    year = {2024},
    volume = {93},
    number = {2},
    pages = {129-143},
}

@article{chousionis_erdogan_tzirakis_2015,
    author = {Chousionis, Vasilis and Erdoğan, M. Burak and Tzirakis, Nikolaos},
    title = {{Fractal Solutions of Linear and Nonlinear Dispersive Partial Differential Equations}},
    journal = {Proceedings of the London Mathematical Society},
    year = {2015},
    volume = {110},
    number = {3},
    pages = {543-564},
}

@book{erdogan_tzirakis_2016,
    author = {Erdo{\u{g}}an, M. Burak and Tzirakis, Nikolaos},
    title = {Dispersive {P}artial {D}ifferential {E}quations: {W}ellposedness and {A}pplications},
    series = {London Mathematical Society Student Texts},
    year = {2016},
    publisher = {Cambridge University Press},
}

\end{document}